\documentclass[11pt]{article}

\usepackage{amsmath, amsfonts, amssymb,latexsym,wasysym,graphicx,stmaryrd,tikz,txfonts}
\input epsf.sty

\newtheorem{Theorem}{Theorem}[section]
\newtheorem{Lemma}{Lemma}[section]
\newtheorem{Proposition}[Lemma]{Proposition}

\newtheorem{Definition}[Lemma]{Definition}
\newtheorem{Remark}[Lemma]{Remark}

\newcommand{\BEQ}{\begin{equation}}     
\newcommand{\BEA}{\begin{eqnarray}}
\newcommand{\BD}{\begin{displaymath}}
\newcommand{\EEQ}{\end{equation}}       
\newcommand{\EEA}{\end{eqnarray}}
\newcommand{\ED}{\end{displaymath}}

\newcommand{\Del}{\Delta}
\newcommand{\eps}{\varepsilon}          

\newcommand{\R}{\mathbb{R}}

\newcommand{\N}{\mathbb{N}}

\def\proba{{\mathbb{P}}}

\def\esper{{\mathbb{E}}}

\newcommand{\eop}{\hfill $\Box$}        

\newcommand{\half}{{1\over 2}}          

\catcode`\@=11
\def\numberbysection{\@addtoreset{equation}{section}
        \def\theequation{\thesection.\arabic{equation}}}
\numberbysection

\begin{document}

\vspace*{1.5cm}
\begin{center}
{\Large \bf Global existence and smoothness for  solutions of viscous Burgers equation. (2) The
unbounded case: a characteristic flow study}

\end{center}

\vspace{2mm}
\begin{center}
{\bf  J\'er\'emie Unterberger$^a$}
\end{center}

\vskip 0.5 cm
\centerline {$^a$Institut Elie Cartan,\footnote{Laboratoire 
associ\'e au CNRS UMR 7502} Universit\'e de Lorraine,} 
\centerline{ B.P. 239, 
F -- 54506 Vand{\oe}uvre-l\`es-Nancy Cedex, France}
\centerline{jeremie.unterberger@univ-lorraine.fr}

\vspace{2mm}
\begin{quote}

\renewcommand{\baselinestretch}{1.0}
\footnotesize
{We show that the homogeneous viscous Burgers equation $(\partial_t-\eta\Del) u(t,x)+(u\cdot\nabla)u(t,x)=0,\ (t,x)\in\R_+\times\R^d$ $(d\ge 1, \eta>0)$
  has a globally defined
smooth solution  if the initial condition $u_0$ is a smooth function  growing  like $o(|x|)$  at infinity. The proof relies mostly on estimates of the random characteristic flow defined by a Feynman-Kac representation of the solution.  Viscosity independent a priori
bounds for the solution are derived from these. The regularity of the solution is then proved for fixed $\eta>0$ using Schauder estimates.

The result extends with few modifications to initial conditions growing abnormally large in regions with small relative volume, separated by well-behaved bulk regions, provided
these are stable under the characteristic flow with high probability. We provide a large
family of examples for which this loose criterion may be verified by hand.
}
\end{quote}

\vspace{4mm}
\noindent


\smallskip


\newpage

\tableofcontents



\section{Introduction and summary of results}


\subsection{Introduction}


The $(1+d)$-dimensional viscous Burgers equation is the following non-linear PDE, 
\BEQ (\partial_t -\eta\Del+u\cdot\nabla )u=0, \qquad u\big|_{t=0}=u_0 \label{eq:Burgers}\EEQ
for a velocity $u=u(t,x)\in\R^d$ ($d\ge 1$), $(t,x)\in\R_+\times\R^d$,
where $\eta>0$ is a viscosity coefficient, $\Del$ the standard Laplacian on $\R^d$,   $u\cdot\nabla u=\sum_{i=1}^d
u_i\partial_{x_i}u$ the convection term, and $g$ a  continuous forcing term.  Among other things, this
fluid equation describes the hydrodynamical limit of interacting particle systems \cite{Spo,KipLan}, is a
simplified version without pression of the incompressible Navier-Stokes equation, and also (adding a random forcing term in the right-hand side) an interesting toy model for the study of turbulence \cite{BecKha}.

The traditional strategy to show a priori estimates for this equation, see e.g. \cite{KL}, 
is to combine integral $L^2$-estimates (the simplest of which coming from the energy
balance equation) with the maximum principle. The latter, valid for any transport
equation -- but not for the related Navier-Stokes equation -- implies a uniform bound for
the supremum $||u_t||_{\infty}$ of the solution,  $||u_t||_{\infty}\le 
||u_0||_{\infty}$.

In a previous article \cite{Unt-Bur1}, we showed that the maximum principle alone was enough to show global existence and boundedness of the solution, provided the initial
solution is {\em bounded} together with its derivatives to order 2. In particular, it is
not necessary to assume that $u_0$ or $g$ are in $L^2$-spaces to solve the equation.
Also, our bounds do not grow exponentially in time, contrary to the classical bounds
based on energy estimates, see e.g.  \cite{KL}.

In the present work, we aim at {\em relaxing the boundedness hypothesis as much as possible}.
If the initial condition is unbounded, then  the maximum principle does not make sense any more.  For
 solutions of some {\em scalar} parabolic equations, e.g. of viscous Hamilton-Jacobi equations,
the comparison principle allows one to define viscosity solutions growing at infinity
\cite{DaLioLey}.  However, here $u$ is not scalar, nor can it be reduced in general to the solution of a Hamilton-Jacobi
equation (save in dimension 1), so it is not at all clear if such
a strategy can work. Instead we tackle the problem {\em from a dynamical system perspective}
and ask ourselves: {\em can one find general criteria ensuring that characteristics of the flow
do not blow up ?} 

It turns out that this question is really the crux of the problem. Let us explain roughly why in the case of zero viscosity ($\eta=0$). Recasting this Eulerian fluid equation into a Lagrangian language,  $u$ is constant along its (time-reversed) characteristics,
 defined as the solutions of the ordinary differential equations
$\frac{d}{ds} x(t;s,x)=u(t-s,x(t;s,x))$ with initial condition $x$; in other words,
$u(t,x)=u(t-s,x(t;s,x))$. In particular $u(t,x)=u_0(x(t;t,x))$ is a priori well defined if
$u_0$ is, no matter how large $u_0$ can be. The argument is clearly faulty as the
characteristic $x(t;s,x)$ may indeed blow up if $u_0$ grows too fast at infinity. This
is clear if one replaces $u$ by the approximation $\tilde{u}$ (denoted $u^{(1)}$ later
on) defined by: $\tilde{u}(t,x):=u_0(\tilde{x}(t;t,x))$,  $\tilde{x}(t;\cdot,x)$
solving the above differential equation, but with the velocity $u(t-s,\cdot)$ 
approximated by the initial velocity $u_0(\cdot)$, namely,  $\frac{d}{ds} \tilde{x}(t;s,x)=u_0(\tilde{x}(t;s,x))$.  This equation does not blow up in finite time if 
$u_0$ is  Lipschitz and has sublinear growth at infinity.  Since linear growth is
really a border case, we shall rather consider as {\em prototypical initial velocity} a function with {\em strictly
sublinear} growth, namely, $|u_0(x)|=O_{|x|\to\infty}(|x|^{1/\kappa})$, $\kappa>1$, for
which $\tilde{x}(t;t,x)$ grows for large time like $t^{\kappa/(\kappa-1)}$.  
But then one may go one step further and remark that the {\em instantaneous value} of
$u_0$ at some point is  not so  important.  Indeed, in one dimension, the non blow-up criterium states that the time needed to go from $x$ to $x'$
(equal to  $\int_{x}^{x'} \frac{dy}{u_0(y)}$ if e.g. $x<x'$ and $u_0>0$) must diverge
when  $|x'|\to\infty$; this does not prevent $u_0$ from becoming arbitrary large
in regions with small relative size, provided these are separated by large bulk intervals
where $u_0$ grows sublinearly and which therefore take up a large time to cross. In short,
we are happy if $\tilde{x}(t;t,x)-x=O(t^{\kappa/(\kappa-1)})$  for $t$ {\em large}.

Surely enough, this last criterion should not be taken seriously for a number of 
obvious reasons (it is dimension-dependent, what $t$ large means is not clear, the
connection to the original non-linear equation is not clear, what happens in case of
non-zero viscosity, etc.), but it really is the inspiration of the present work. 
Let us sketch the  answer to  some of the objections we have just raised.
First, as in \cite{Unt-Bur1}, we use the following scheme of successive approximations to
the solution. We solve inductively the 
linear transport equations,
\BEQ u^{(-1)}:=0; \EEQ
\BEQ (\partial_t-\Del+u^{(m-1)}\cdot\nabla)u^{(m)}=0, \ \ u^{(m)}\big|_{t=0}=u_0 \qquad (m\ge 0). \label{eq:um}\EEQ
If the sequence $(u^{(m)})_m$ converges locally in $C^{1,2}$-norms, then the limit is a fixed point of (\ref{eq:um}),
hence solves the Burgers equation.
The Feynman-Kac formula implies the following well-known representation of the solution
of (\ref{eq:um}) in terms of {\em random characteristics} $X^{(m)}(t,\cdot)$,
\BEQ u^{(m)}(t,x)=\esper[u_0(X^{(m)}(t,x))], \label{eq:FK1} \EEQ
where $X^{(m)}(t,x):=X^m(t;t,x)$ is the solution at time t of a stochastic differential
equation driven by a standard Brownian motion $B$,
\BEQ dX^{(m)}(t;s,x)=u^{(m-1)}(t-s,X^{(m)}(t;s,x))ds+dB_s,  \label{eq:FK2} \EEQ
 started at $X^{(m)}(t;0,x):=x$.

In section 2 we concentrate on {\em prototypical initial velocities}, i.e. study
Burgers equation under the hypothesis 
\BEQ  {\mathrm{(Hyp1)}}\qquad |u_0(x)|\le U (1+|x|)^{1/\kappa}, \qquad x\in\R^d
\label{intro:Hyp1} \EEQ
with $\kappa>1$ and $U\ge 1$. Solving for the random
characteristic $X^{(1)}$ (which coincides with the above deterministic characteristics  $\tilde{x}$ in the zero viscosity case), we prove that for $t$ {\em large}, {\em with high probability},
 \BEQ |X^{(1)}(t;s,x)-x|=O\left( \max\left( (Ut)^{\kappa/(\kappa-1)}, Ut |x|^{1/\kappa}
 \right)\right), \label{intro:X1} \EEQ
  thus   retrieving for $t$ large  the behaviour in $O(t^{\kappa/(\kappa-1)})$). Then we note that $X^{(m)}, m\ge 2 $ solves essentially the same equation as $X^{(1)}$
 since $u^{(m-1)}(t-s,y)=\tilde{\esper}\left[ u_0(X^{(m-1)}(t-s,y)) \right]$
 is the average of $u_0$ on some weighted cloud of points in a neighbourhood of $y$.
At this point it is natural to introduce what we call a  {\em generalized flow with initial velocity}
$u_0$ (see Definition \ref{def:generalized-flow}). Roughly speaking, at least in the non-viscous case, this is an ordinary differential equation of the form
$\frac{d}{ds}y(t;s,x)=u_0({\cal X}(t;s,y(t;s,x)))$ where ${\cal X}(t;s,\cdot)$ satisfies an
estimate of the same form as $X^{(1)}(t;s,\cdot)$ (see eq. (\ref{intro:X1})). In the
{\em viscous case}, we first convert the {\em stochastic} differential equation (\ref{eq:FK2}) into an {\em ordinary} differential equation {\em with random coefficients}
by subtracting the additive noise $B$ (see section 2.3). Then  {\em viscous generalized
flows}  (see Definition \ref{def:viscous-generalized-flow}) are (non-viscous) generalized flows, in which spatial arguments have been translated by the noise. Now the interesting
property about  generalized flows $y(t;\cdot,x)$ is that they themselves satisfy 
some version of (\ref{intro:X1}), where $U$ is the constant appearing in (Hyp1)
 (see Lemmas \ref{lem:generalized-flow-bound}, \ref{lem:viscous-generalized-flow-bound}) . As a result,  we are able to obtain inductively bounds for $X^{(m)}$ of the type (\ref{intro:X1}) which are {\em uniform in} $m$.

At this point, one would be tempted to define an {\em admissible} initial velocity as
a function $u_0$ for which the inductive Lemmas \ref{lem:generalized-flow-bound}, \ref{lem:viscous-generalized-flow-bound} hold.  {\em As pointed out above, the restriction
'for $t$ large' is essential}: should we require that (\ref{intro:X1}) hold for $t$ small,  this
would directly imply a sublinear bound on the velocity. Actually,  working out the computations, it appears very soon that $t\gtrsim U^{-1}$ is the right condition.
Now, while for a given function $u_0$  the conclusions of Lemmas \ref{lem:generalized-flow-bound}, \ref{lem:viscous-generalized-flow-bound} may be eventually verified by hand,
it turns out that, leaving aside the settled case of functions satisfying (Hyp1), 
it is difficult to produce any interesting example of admissible velocity. The
reason is of topological origin: we need some criterion ensuring inductively the
{\em  stability} under the characteristic flows of the {\em safe zones} where $u_0$ is sublinear. To be more specific (see section 3), we {\em assume} that $u_0$ is sublinear in some 'bulk'
{\em safe} region $\cal S$ (connected or not), while it is essentially arbitrary
in a countable disjoint union of 'thin' {\em dangerous} regions $({\cal A}_i)_{i\in I}$. 
In Definition \ref{def:dangerous-zones} we choose these to be annuli, but clearly
this is only a reasonable, practical choice. The important thing is that, sticking to the non-viscous case for the time being, {\em provided} the safe zones are 'fat' enough, one is  able to prove inductively a {\em safe zone stability property} stating that 
$$\left( x^{(m-1)}(t;s,x) \in {\cal S}(t-s), t\ge s\ge 0 \right)\Longrightarrow \left( x^{(m)}(t;s,x) \in {\cal S}(t-s), t\ge s \ge  0 \right),$$ 
where $t\mapsto {\cal S}(t)$ is some decreasing family of non-empty subsets  with
${\cal S}(0)={\cal S}$ (see
Theorem \ref{th:admissible}).  In this way we show that  $x^{(m)}(t;s,x)\in {\cal S}(t-s)$ for all $m$ as soon as $x\in {\cal S}(t)$. Let ${\cal A}(t):=\R^d\setminus {\cal S}(t)$ be the {\em enlarged dangerous zone}. If $x\in {\cal A}(t)$, then $x$  may a priori jump to the boundary of ${\cal A}(t)$ in arbitrarily short time, after which it cannot
escape from the safe zone any more due to the safe zone stability property. {\em If} ${\cal A}(t)$ is
still a disjoint union $({\cal A}_i(t))_{i\in I}$ of thin regions, then this may (and does under our assumptions for $({\cal A}_i)_{i\in I}$) prove enough to show a uniform
bound of the type (\ref{intro:X1}). Thus the safe zone stability property is an efficient
replacement for the inductive property of Lemmas \ref{lem:generalized-flow-bound}, \ref{lem:viscous-generalized-flow-bound}.

A straigthforward generalization of these arguments to the {\em viscous case} appears
to be impossible at first sight, since one may always fall into the dangerous zone by translating by some random amount the spatial
arguments. Even though these random amounts are bounded {\em in average}, without
additional assumptions on $u_0$, it may happen, with a small but nonzero probability,
that random characteristics blow up. So much for the debit side. On the credit side,
one sees that  the translation by  random paths $(B_t)_{t\ge 0}$ bounded by $o(t)$ for $t$ large (which is the case of the overwhelming majority of random paths since $B_t$ is roughly
of order $\sqrt{t}$) should not affect the usual displacement bound in
$ O\left( \max\left( (Ut)^{\kappa/(\kappa-1)}, Ut |x|^{1/\kappa}
 \right)\right)$ (see eq. (\ref{intro:X1})) since $o(t)\ll (Ut)^{\kappa/(\kappa-1)}$
 for $t\ge U^{-1}$ (see (\ref{eq:U1/2}) for a more quantitative statement). In short,
 as emphasized in section 2.3, {\em convection prevails over diffusion} in normal
 conditions. Since the opposite regime where diffusion prevails over convection is
 highly improbable, only very mild assumptions (e.g. polynomial growth at infinity) are $u_0\big|_{{\cal A}}$ is required
to extend the safe zone stability property argument to the viscous case. A precise
statement may be found e.g. in Lemma \ref{lem:viscous-admissible}.

\medskip

Once one has a uniform control of the random characteristics, and some polynomial
'a priori' bound on $u_0$, one may start about proving the convergence of the scheme
(\ref{eq:um}), which is the subject of section 4. From that point on, we
follow a more conventional  course of action, which  is
sketched in the next paragraph.


\subsection{Summary of results}


The general assumptions on the initial velocity $u_0$ are written down in the preamble
of section 4. Fix $U\ge 1, \kappa>1$. We demand the following: (i) $u_0$ is  $C^2$;  (ii) $u_0$, $\nabla u_0$ and
$\nabla^2 u_0$   grow at most polynomially at infinity (these we call {\em a priori
bounds} for $u_0$, see (\ref{eq:a-priori-bound})); plus a third condition (iii)
stating roughly that the characteristic flows $s\mapsto X^{(m)}(t;s,x)$ may be estimated for
$t\gtrsim U^{-1}$ like
the deterministic flow $s\mapsto y(s,x)$ defined by the ordinary differential equation
$\frac{d}{ds}y(s,x)=(1+|y(s,x)|)^{1/\kappa}$  with initial condition $y(0,x)=x$, 
except when $\sup_{0\le s\le t} |B_s|$ overrides the usual displacement bound (\ref{intro:X1}), the latter condition defining the so-called highly improbable {\em abnormal regime} where
diffusion prevails over convection. Depending on whether one wants {\em examples} built following the
above arguments (with explicit 'safe' and 'dangerous' zones, etc.)  which are {\em 
sufficient} to ensure such estimates, or one rather looks for more or less 'necessary' conditions {\em a minima}
on the characteristics in the abnormal regime ensuring that all subsequent estimates
(on $u^{(m)},\nabla u^{(m)}$...) remain unaffected, one obtains different versions of
(iii). The {\em sufficient} condition (iii) is based on Definition \ref{def:dangerous-zones}:

\medskip

\noindent {\bf Theorem 1} {\em (see Definition \ref{def:dangerous-zones}, Theorem \ref{th:viscous-admissible} and (\ref{eq:bound-u}))
Let $(R_n)_{n\ge 1}$ be an increasing sequence, $1\le R_1<R_2<R_3<\ldots$ such that, for all $i\ge 1$, 
\BEQ R_{2i}-R_{2i-1}\le R_{2i-1}^{1/\kappa}, \label{intro:R1} \EEQ
\BEQ  R_{2i+1}\ge 4 R_{2i} \label{intro:R2}
.\EEQ

Let $\tilde{u}_0:\R^d\to\R^d$ be an initial velocity satisfying {\em (Hyp1)} (see
(\ref{intro:Hyp1})) for some constants
$U\ge 1$, $\kappa>1$. Let $u_0:\R^d\to\R^d$ be any Lipschitz function coinciding with
$\tilde{u}_0$ outside the union of annular 'dangerous zones' $\cup_{i\ge 1} {\cal A}_i$, ${\cal A}_i:=B(0,R_{2i})\setminus B(0,R_{2i-1})$, and satisfying the a priori bounds
(\ref{eq:a-priori-bound}). Let also $M_t:=1+ \frac{\sup_{0\le s\le t} |B_s|}{\sqrt{t}}$. Then  the sequence of noise-translated characteristics
$(Y^{(m)}(t;\cdot,x))_{m\ge 0}$, $Y^{(m)}(t;s,x):=X^{(m)}(t;s,x)-B_s$, satisfies the following uniform in $m$ estimates:

\BEQ |Y^{(m)}(t;s,x)-x|\lesssim   \langle Ut\rangle  \max(\langle Ut\rangle^{\kappa/(\kappa-1)},|x|)^{1/\kappa} \qquad  {\mathrm{if}} 
\ M_t\sqrt{t}\le \max(\langle Ut\rangle ^{\kappa/(\kappa-1)}, \langle Ut \rangle \langle x\rangle^{1/\kappa}); \label{intro:kappa0} \EEQ

in the {\em normal regime}, otherwise 
\BEQ  |Y^{(m)}(t;s,x)-x|\lesssim \left(\frac{M_t\sqrt{t}}{\langle Ut\rangle}\right)^{\kappa} \label{intro:kappa} \EEQ

Furthermore, estimates (\ref{intro:kappa0}), (\ref{intro:kappa}) imply for
$u^{(m)}$, $m\ge 0$ defined by Feynman-Kac's formula (\ref{eq:FK1})

\BEQ |u^{(m)}(t,x)| \lesssim  K_0 (|x|+\langle Ut\rangle^{\kappa/(\kappa-1)})^{\frac{\alpha}{2}+\frac{1}{\kappa}}. \label{intro:bound-u} \EEQ
}

On the other hand, bounds (\ref{intro:kappa}) in the {\em abnormal regime}
$M_t\sqrt{t}\ge \max(\langle Ut\rangle ^{\kappa/(\kappa-1)}, \langle Ut \rangle \langle x\rangle^{1/\kappa})$ may be considerably softened without harming ulterior bounds.
In particular, substituting to (\ref{intro:kappa}) the condition
\BEQ  |Y^{(m)}(t;s,x)-x|\lesssim (M_t\sqrt{t})^{\kappa'} \label{intro:kappa'} \EEQ
for some arbitrary exponent $\kappa'\ge 1$, one still has (\ref{intro:bound-u}). 
Demanding only (\ref{intro:kappa0}) and (\ref{intro:kappa'}), we get our 'necessary' condition (iii').  Of course, it remains to be proved that there are different choices of dangerous zones -- or, from a wider perspective, of 
functions $u_0$ -- for which (\ref{intro:kappa'}) holds but not (\ref{intro:kappa}).
In any case, bounds in section 4 are based on (\ref{intro:kappa'}).

\medskip

Let us comment on conditions (\ref{intro:R1}), (\ref{intro:R2}). Condition (\ref{intro:R1}) states that the {\em width of the dangerous zone} ${\cal A}_i$ is
{\em smaller} than the expected displacement $0\left(\max\left( (Ut)^{\kappa/(\kappa-1)}, Ut |x|^{1/\kappa}
 \right) \right)$ (see (\ref{intro:X1})) for all $t\ge U^{-1}$. Condition 
 (\ref{intro:R2}) states that the {\em width of the safe zone} $B(0,R_{2i+1})\setminus
 B(0,R_{2i})$ is {\em larger} than the expected displacement for $|x|\gg \langle Ut\rangle^{\kappa/(\kappa-1)}$. The latter condition (characteristic of the so-called
 {\em short-time regime}, where $\max\left( (Ut)^{\kappa/(\kappa-1)}, Ut |x|^{1/\kappa}
 \right) \lesssim |x|$)  comes up naturally right from
 the beginning (see section 2.1). There is nothing special about the coefficient 4
 in (\ref{intro:R2}), and our results carry through if $R_{2i}-R_{2i-1}\le CR_{2i-1}^{1/\kappa}, R_{2i+1}\ge (1+\eps)R_{2i}$ with $C,\eps>0$ arbitrary, but then implicit
 constants also depend on $C,\eps$, instead of depending only on the dimension $d$ and
 on the exponents $\kappa,\kappa'$.

 \medskip

From a logical point of view, the above  Theorem  is inaccurate since it provides
a priori bounds for objects such as $Y^{(m)}(\cdot;\cdot,\cdot)$, $u^{(m)}(\cdot,\cdot)$
without proving their existence. In particular, one must prove inductively that 
$(u^{(m)})_{m\ge 0}$ are $C^1$, so that the transport equations (\ref{eq:um}) are 
well-posed and we can use Cauchy-Lipschitz's theorem to define uniquely the characteristics. Ultimately we prove the following:

\medskip

\noindent {\bf Theorem 2} (see  sections 4.2, 4.3, 4.4)
 {\em Assume that hypotheses  (\ref{intro:kappa0}), and (\ref{intro:kappa}) (or more generally (\ref{intro:kappa'})) hold, and that $u_0,
 \nabla u_0,\nabla^2 u_0$ satisfy the following a priori bounds (see (\ref{eq:a-priori-bound})),
 \BEQ |u_0(x)|\le K_0(1+|x|)^{\frac{\alpha}{2}+\frac{1}{\kappa}}, 
\qquad |\nabla u_0(x)|\le  K_1 (1+|x|)^{\alpha+\frac{2}{\kappa}}, \qquad 
|\nabla^2 u_0(x)|\le K_2 (1+|x|)^{\frac{3}{2}(\frac{\alpha}{2}+\frac{1}{\kappa})} 
\label{intro:a-priori-bound} \EEQ
with 
\BEQ K_0\le U^{\frac{\beta}{2}+1}, \qquad K_0\le K_1^{1/2},  \qquad U\le K_1\le K_2^{2/3}, \label{intro:K} \EEQ
 for
some exponents $\alpha,\beta\ge 0$.
 
Let
$v^{(0)}:=u^{(0)}$ and $v^{(m)}:=u^{(m)}-u^{(m-1)}$ ($m\ge 1$). Fix $\gamma\in(0,1)$. Then there exists
a universal constant $C=C(d,\kappa,\kappa',\alpha,\beta,\gamma)>1$ such that, for $m\ge 0$,

\BEQ |\nabla u^{(m)}(t,x)|\le C^2 K_1 (|x|+ \langle Ut\rangle^{\kappa/(\kappa-1)})^{\alpha+\frac{2}{\kappa}};  \label{intro:nabla-u} \EEQ
\BEQ |\nabla^2 u^{(m)}(t,x)|\le C^4 K_2 (|x|+ \langle Ut\rangle^{\kappa/(\kappa-1)})^{3(\frac{\alpha}{2}+\frac{1}{\kappa})}; \label{intro:nabla2-u}  \EEQ
\BEQ |v^{(m)}(t,x)|\le CK_0 (t/mT_{min}(t,x))^m (|x|+ \langle Ut\rangle^{\kappa/(\kappa-1)})^{\alpha+\frac{2}{\kappa}};  \label{intro:v} \EEQ
\BEQ |\nabla v^{(m)}(t,x)|\le C^3 K_2^{2/3} (t/m\tilde{T}_{min}(t,x))^{\gamma m/2} (|x|+ \langle Ut\rangle^{\kappa/(\kappa-1)})^{\alpha+\frac{2}{\kappa}}  \label{intro:nabla-v} \EEQ
where 
\BEQ T_{min}(t,x):=\left(C^3 K_1  (|x|+ \langle Ut\rangle^{\kappa/(\kappa-1)})^{\alpha+\frac{2}{\kappa}} \right)^{-1}, \tilde{T}_{min}(t,x):=\left(C^3 K_2^{2/3}  (|x|+ \langle Ut\rangle^{\kappa/(\kappa-1)})^{\alpha+\frac{2}{\kappa}} \right)^{-1}.\EEQ
}

\medskip\noindent  Estimates (\ref{intro:v}, \ref{intro:nabla-v}) imply convergence in absolute value
of the series $\sum_{m\ge 0} v^{(m)}$, $\sum_{m\ge 0} \nabla v^{(m)}$, from which
it may be concluded by standard arguments that the limit $v$ satisfies Burgers' equation.
Theorems 1 and 2 must actually be proved simultaneously since they are based on
induction (the a priori bounds at rank $m-1,m$ proved in Theorem 1 are used to
prove rank $m$ gradient estimates (\ref{intro:nabla-u}) of Theorem 2, from which one can justify the a priori bounds at rank $m+1$, etc.)

\medskip

Let us comment on a priori bounds (\ref{intro:a-priori-bound}), and in particular on
(\ref{intro:K}). As noted in our previous article \cite{Unt-Bur1}, dimensional analysis,
confirmed by the initial perturbative expansion but also by Schauder estimates for large $t$,    tells us that $u$, $\nabla u$,
$\nabla^2$ should scale like $L^{-1}$, $L^{-2}$, $L^{-3}$ for some reference length $L$ depending
on the initial condition, at least for bounded solutions. (In our setting where $u_0$
may increase polynomially, we have included an extra reference length $\approx 1$.)
This account for the relations between the exponents appearing in (\ref{intro:a-priori-bound}), (\ref{intro:K})), except for $\beta$ which is arbitrary. Note that  $\beta$ does not appear in the bounds (\ref{intro:nabla-u},\ref{intro:nabla2-u},\ref{intro:v},\ref{intro:nabla-v}), except in the numerical constant $C$. Finally the hypotheses
$K_0\le K_1^{1/2}, U\le K_1\le K_2^{2/3}$ may be discarded provided one defines as in \cite{Unt-Bur1} some
 constant $K:=\max(U,K_0^2,K_1,K_2^{2/3})$ homogeneous to an inverse length, and replaces $K_0,K_1,K_2$ in
(\ref{intro:nabla-u},\ref{intro:nabla2-u},\ref{intro:v},\ref{intro:nabla-v}) by
$K^{1/2},K,K^{3/2}$, thus equating $\tilde{T}_{min}$ with $T_{min}$. 
    
\medskip

Let us finally say some words about the strategy of proof (see section 4.1 for more details), 
which follows closely that of our previous article \cite{Unt-Bur1}. In principle,
we would like to prove  the gradient bounds (\ref{intro:nabla-u}), (\ref{intro:nabla2-u}), (\ref{intro:nabla-v}) by using Feynman-Kac's formula and hypotheses (\ref{eq:a-priori-bound}), (\ref{intro:kappa0}), (\ref{intro:kappa'}) in an initial regime $t\le 
T_{min}(0,x)= (C^3 K_1 (1+|x|)^{\alpha+\frac{2}{\kappa}})^{-1}$, beyond which 
exponential factors due to separation of trajectories become large. However this makes
no sense in itself  since $T_{min}(0,x)\to_{|x|\to\infty} 0$. Furthermore, we are not even able to prove such estimates  if one takes into account the contribution of the 'abnormal
regime' to the expectation appearing in Feynman-Kac's formula. The solution to these
problems is to rewrite $u^{(m)}$ as the sum of a series with general term $\underline{u}^{(m,n)}:=u^{(m,n)}-u^{(m,n-1)}$, where $u^{(m,n)}$, $n\ge 0$ solves a {\em penalized} transport equation meant as a
smoothened substitute of the original equation solved on the dyadic ball $B(0,2^n)$
 (see section 4.2). Then $\nabla u^{(m,n)}$, and similarly $\nabla^2 u^{(m,n)},\nabla v^{(m)}$ may be proved inductively to satisfy
 (\ref{intro:nabla-u},\ref{intro:nabla2-u},\ref{intro:nabla-v}) for $t\le T_n:=(C^3 K_1 (2^n)^{\alpha+\frac{2}{\kappa}})^{-1}
 \approx T_{min}(0,2^n)$.  Furthermore, for $x$ {\em small}, namely, if $|x|\ll 2^n$,
 then Gaussian bounds for Brownian motion imply that $\nabla u^{(m,n)}(t,x),\nabla^2 u^{(m,n)}(t,x),\nabla v^{(m)}$ are exponentially small; intuitively this is clear since
the only contribution to $\nabla\underline{u}^{(m,n)}$ comes from characteristics
$X^{(m)}(t;\cdot,x)$ which go very far away from $x$, crossing the boundary of $B(0,2^n)$.
 Extension of these bounds to larger $t$ is proved using
 home-made (interior) Schauder estimates proved in our previous article \cite{Unt-Bur1}.

 Finally, the series in $n$ converge thanks to the estimates in the small $x$ regime.

\bigskip

{\em Notations:} we let $\langle t\rangle:=\max(1,t)$ for $t\in\R_+$, $\langle x\rangle:=
\max(1,|x|)$ for $x\in\R^d$. Also, given two functions $f,g$, $f\lesssim g$ (resp.
$f\gtrsim g$) means: there exists an overall constant $C$ (depending only on $d$  and on 
the exponents $\kappa,\kappa',\alpha,\beta,\gamma$
possibly) such that $|f(x)|\le C |g(x)|$ (resp. $|f(x)|\ge C|g(x)|$) on the set where $f,g$ are defined. Then $f\approx g$ means: $f\lesssim g$ and $f\gtrsim g$. 


\section{A prototypical example}


In this section we are only interested in providing a priori bounds for the
random paths $X^{(m)}(.;.,.)$, {\em assuming} that the sequence of transport
equations (\ref{eq:um}) admits a unique smooth solution represented by
Feynman-Kac's formula (\ref{eq:FK1},\ref{eq:FK2}). By rescaling we assume $\eta=1$
({\em viscous case}) or $\eta=0$ ({\em non-viscous case}), the latter case serving
essentially as an illustration.

{\em \noindent We assume throughout that $u_0$ is $C^1$}; this is a priori not absolutely necessary (because of the regularizing properties of the heat kernel), but reasonable if one wants to define properly the random characteristics down to time $0$. We make here the following hypothesis:

\medskip

{\centerline{\bf (Hyp1) There exist constants $U\ge 1$, $\kappa>1$ such that  \qquad $|u_0(x)|\le U(1+|x|)^{1/\kappa}$.}}

\medskip

The condition $U\ge 1$ is of course inessential; it avoids having to distinguish
between the factors $O(U)$ and the factors $O(1+U)$ which pop up in the proofs. 
Assuming $u_0$ is small, optimal results using our arguments may be obtained by rescaling the solution and the
time-variable in such a way that $\sup_{x\in\R^d} \frac{|u_0(x)|}{(1+|x|)^{1/\kappa}}=1$, but
mind that this reintroduces a viscosity parameter into the story,  producing in turn a time
rescaling in the bounds (which is very easy to write down by following the
computations step by step).  

\medskip

A prototypical family of natural examples is of course smooth functions $u_0$
satisfying $u_0(x)=F(\frac{x}{|x|})U|x|^{1/\kappa}$ outside $B(0,1):=\{x\in\R^d \ |\ |x|<1\}$,
where $F:S^d\to S^d$ is a smooth function preserving the sphere $S^d:=\{
|x|=1\}$. 

\medskip

In section 3 we shall see that a priori bounds similar to those shown in this section
may be obtained for much more general initial data.


\subsection{Generalities}


We study in this paragraph the flows of ordinary differential equations ({\em ode}'s
for short) of the type
$\dot{x}=u_0(x)$ where $u_0$ satisfies (Hyp1) with parameters $U,\kappa$ such that $U\ge 1$,
 $\kappa>1$.

 We start by introducing a family of typical  ode's depending on a parameter $x_{min}\ge 0$ which we call {\em cut-off}.  

\begin{Definition} $(x_{min}>0)$
Let $\Phi_{\kappa,U,x_{min}}(t,x)$ be the solution at time $t\ge 0$ of the scalar ode 
$\frac{d}{dt}x(t)=U(x_{min}+|x(t)|)^{1/\kappa}$ started at $x(0)=x\in\R$.
\end{Definition}

Solving for $x\ge 0$, one gets ($\kappa>1$ is of course necessary to get a global solution)
\BEQ x(t)=\left( (x+x_{min})^{\frac{\kappa-1}{\kappa}}+\frac{\kappa-1}{\kappa} Ut \right)^{\frac{\kappa}{\kappa-1}} - x_{min}, \qquad t,x\ge 0 \label{eq:phikappaU}.\EEQ
The above solution extends to $t\le 0$ or $x\le 0$ as follows.
If $x\le 0$, $\Phi_{\kappa,U,x_{min}}(t,x)$ reaches $0$ after a time 
 $t=T_{\kappa,U}(x)=
U^{-1} \frac{\kappa}{\kappa-1} \left( (x_{min}+|x|)^{(\kappa-1)/\kappa}
-x_{min}^{(\kappa-1)/\kappa} \right) $, after which we define
$\Phi_{\kappa,U}(t,x):=\Phi_{\kappa,U}(t-T_{\kappa,U}(x),0)>0$. Then (by symmetry)
$\Phi_{\kappa,U}(-t,-x)=-\Phi_{\kappa,U}(t,x)$.

By convention, we let $\Phi_{\kappa,U}(t,x)=\lim_{x_{min}\to 0^+} \Phi_{\kappa,U,x_{min}}(t,x).$

\medskip

The ode's we are interested in are ode's on $\R^d$. Fix $U\ge 1$ and $\kappa>1$.

\begin{Definition} \label{def:velocity}
An ode $\frac{d}{ds}x(s)=v(s,x(s))$ in $\R^d$ has {\em velocity bounded by
$U(x_{min}+|\cdot|)^{1/\kappa}$} on $[0,t]$ if\\ $|v(s,y)|\le U(x_{min}+|y|)^{1/\kappa}$ for all 
$s\in[0,t]$ and $y\in\R^d$. 

\noindent If the velocity field $v$ satisfies this property, we write $v\in {\cal V}_{\kappa,U,x_{min}}(t)$.
\end{Definition}

\begin{Definition} \label{def:BkappaU} Let $B_{\kappa,U,x_{min}}(t,x):=\cup_{v\in {\cal V}_{\kappa,U,x_{min}}(t)} \Big\{ (x(s))_{0\le s\le t} \ |\ (x(s))_{0\le s\le t}$ solution\\  of
the ode $\frac{d}{ds}x(s)=v(s,x(s))$ started at $x(0)=x\Big\}$. 
Let also
\BEQ B_{\kappa,U}(t,x):=\cup\{B_{\kappa,U,x_{min}}(t,x); x_{min}\le 1\}.\EEQ
\end{Definition}

Let us first study $B_{\kappa,U}(t,x)$. 
If the ode $\frac{d}{ds}x(s)=v(s,x(s))$ started at $x$ has a velocity bounded by
$U(1+|\cdot|)^{1/\kappa}$, then  $\frac{d}{ds}| x(s)|\in [-U(1+|x(s)|)^{1/\kappa},
U (1+|x(s)|)^{1/\kappa}]$. Thus $B_{\kappa,U}(t,x)\subset B(x,R_t(|x|))$, where
$R_t(|x|)=\max(x_+(t)-|x|,|x|-x_-(t))$ and $x_{\pm}(t)$ are the solution at time $t$ of the
scalar ode's $\frac{d}{ds}x_+(s)=U(1+|x_+(s)|)^{1/\kappa}$, resp. 
$\frac{d}{ds}x_-(s)=-U(1+|x_-(s)|)^{1/\kappa}$ started at $|x|$.

The reader may easily check by solving either of these  ode's and comparing to
(\ref{eq:phikappaU})  that 
$R_t(|x|)\approx \max(\Phi_{\kappa,U}(t,|x|)-|x|,|x|-\Phi_{\kappa,U}(-t,|x|))$ as soon as
$|x|\gtrsim 1$ or $Ut\gtrsim 1$. Then clearly $|x|-\Phi_{\kappa,U}(-t,|x|)\le \Phi_{\kappa,U}(t,|x|)-|x|$. In absolute generality, it holds $R_t(|x|)\lesssim \Phi_{\kappa,U}\left(\max(t,U^{-1}),|x|\right)-|x|$; the short-time regime $t\lesssim U^{-1}$ is rather uninteresting and need
not be discussed in greater details. Looking more closely at the solution $x(t)$ of 
(\ref{eq:phikappaU}) with $x_{min}\le 1$, we see that there are two regimes, the {\em long-time regime}
where $|x|\ll |Ut|^{\kappa/(\kappa-1)}$ and 
\BEQ |x|\ll U|t|\ |x|^{1/\kappa} \ll |x(t)-x|\approx |x(t)|\approx 
|Ut|^{\kappa/(\kappa-1)}, \EEQ  and the opposite {\em short-time
regime}, $|x|\gg |Ut|^{\kappa/(\kappa-1)}$, where 
\BEQ |Ut|^{\kappa/(\kappa-1)} \ll |x(t)-x|\approx U|t|\  |x|^{1/\kappa}\ll |x| \EEQ
 is small.
Note that 
\BEQ |x(t)-x|\lesssim \max\left( |Ut|^{\kappa/(\kappa-1)},U|t|\, |x|^{1/\kappa}\right) \EEQ
for all values of $t$ and $x$. 

All these estimates generalize straightforwardly to {\em small cut-offs}, $x_{min}\lesssim
(Ut)^{\kappa/(\kappa-1)}$: namely, for such values of $x_{min}$,  $x(s)\in B_{\kappa,O(U)}(t,x)$ for $s\in [0,t]$, as easily shown from the previous computations.

\medskip Things get different when $x_{min}$ is {\em large}, say, $x_{min}>(Ut)^{\kappa/(\kappa-1)}$. Taylor expanding (\ref{eq:phikappaU}) started  from $x>0$, one sees that,
for all $t>0$,
\BEQ x(t)=(x+x_{min})(1+O(Ut\  x_{min}^{-(\kappa-1)/\kappa}))-x_{min}=x+O(Ut\  x_{min}^{1/\kappa}), \qquad x\le x_{min} \EEQ
while
\BEQ x(t)=(x+x_{min})(1+O(Ut\  x^{-(\kappa-1)/\kappa}))-x_{min}=x+O(Ut\  x^{1/\kappa}), \qquad x\ge x_{min} \EEQ
Though we still get two different regimes, it makes sense to say that the long-time
regime has been 'swallowed' by the short-time regime.

Summarizing, we get:

\begin{Lemma} \label{lem:BkappaU}
Let $t\ge 0$ and  $x\in\R^d$.
\begin{itemize}
\item[(1)] ({\em small cut-off regime}) Let  $x_{min}\in[0,(Ut)^{\kappa/(\kappa-1)}]$. Then $B_{\kappa,U,x_{min}}(t,x)\subset B(x,C(\Phi_{\kappa,U}(t,|x|)-|x|))$ for some constant $C\ge 1$. Furthermore, there exists some constant 
$C'\ge 1$ such that, independenly of $x_{min}$ :
\begin{itemize}
\item[(i)] if $|x|\lesssim (Ut)^{\kappa/(\kappa-1)}$ ({\em long-time regime}),
\BEQ B_{\kappa,U,x_{min}}(t,x)\subset B(0,C'(Ut)^{\kappa/(\kappa-1)}); \label{eq:long-time-regime} \EEQ
\item[(ii)] if $|x|\gtrsim (Ut)^{\kappa/(\kappa-1)}$ ({\em short-time regime}),
\BEQ B_{\kappa,U,x_{min}}(t,x)\subset B(x,C'(Ut) |x|^{1/\kappa}). \label{eq:short-time-regime} \EEQ
\end{itemize}

\item[(2)] ({\em large cut-off regime}) There exists some constants $0<c<1<C$ such
that the following holds. Let   $x_{min}\ge (Ut)^{\kappa/(\kappa-1)}$.
Then 
\BEQ B(x,cUt\max(x_{min},|x|)^{1/\kappa})\subset B_{\kappa,U,x_{min}}(t,x)\subset B(x,CUt\max(x_{min},|x|)^{1/\kappa}). \label{eq:large-cut-off-regime} \EEQ
\end{itemize}
\end{Lemma}

Note the following particular case of (\ref{eq:large-cut-off-regime}),
\BEQ B(x,c(Ut)^{\kappa/(\kappa-1)}) \subset B_{\kappa,U,(Ut)^{\kappa/(\kappa-1)}}(t,x)
\subset B(x,C(Ut)^{\kappa/(\kappa-1)}), \qquad |x|\le (Ut)^{\kappa/(\kappa-1)}. 
\label{eq:large-cut-off-regime2} \EEQ

\begin{Remark} \label{rem1}

In particular, an ode with velocity 
\BEQ |v(s,y)|\lesssim
U\left(1+|y|+O(\sqrt{t})\right)^{1/\kappa} \EEQ 
is covered  by Lemma \ref{lem:BkappaU} (1) for $t\ge U^{-1}$ since 
\BEQ \sup_{t\ge U^{-1}} 
\sqrt{t}/(Ut)^{\kappa/(\kappa-1)}=U^{-1/2}\le 1 \label{eq:sqrt-t}. \EEQ
 Perturbation in $O(\sqrt{t})$ do appear as
an effect due to diffusion (see \S 2.3). Thus the general philosophy is that {\em convection prevails over
diffusion} in our setting.
\end{Remark}


\subsection{The non-viscous case}


We set the viscosity $\eta$ to $0$ in this
paragraph. Namely, the zero-viscosity case is interesting in itself, easier to study, and contains already
the main features of the viscous case (see \S 2.3 below). 
We are thus led to consider the approximation sheme 
\BEQ \phi^{(-1)}:=0;\EEQ
\BEQ (\partial_t+\phi^{(m-1)}(t,x)\cdot\nabla)\phi^{(m)}(t,x)=0 , \qquad 
\phi^{(m)}\big|_{t=0}=u_0 \qquad (m\ge 0)
\EEQ
to the non-viscous Burgers equation
\BEQ (\partial_t+\phi\cdot\nabla)\phi=0,  \qquad \phi\big|_{t=0}=u_0 \EEQ
 with initial condition $u_0$ satisfying (Hyp1). The zero-viscosity
Feynman-Kac expression for the solution (compare with (\ref{eq:FK1}), (\ref{eq:FK2})) is given in terms of deterministic characteristics $x^{(m)}(\cdot,x)$, $m\ge 0$,  viz.
\BEQ \phi^{(m)}(t,x)=u_0(x^{(m)}(t,x)),\EEQ
where  $x^{(m)}(t,x):=x^{(m)}(t;t,x)$ is the solution at time $t$ of the ode
\BEA  \frac{d}{ds} x^{(m)}(t;s,x) &=& \phi^{(m-1)}(t-s,x^{(m)}(t;s,x)) \nonumber\\
&=& u_0\big( x^{(m-1)}(t-s, x^{(m)}(t;s,x))\big) \label{eq:xmxm-1} \EEA
with initial condition $x^{(m)}(t;0,x)=x$.
(Later on -- see section 3 -- we shall check inductively that $\phi^{(m)}(t,x)$ is
continuous in time and Lipschitz in $x$, so that (\ref{eq:xmxm-1}) has a unique
solution, possibly only for small time.)

\medskip In particular, 
\BEQ x^{(0)}(t,x)=x;\EEQ
\BEQ \frac{d}{ds} x^{(1)}(t;s,x)=u_0(x^{(1)}(t;s,x)).\EEQ
The ode for $x^{(1)}$ has by (Hyp1) a velocity bounded by $U(1+|\cdot|)^{1/\kappa}$, so, by
Definition \ref{def:BkappaU}, 
\BEQ x^{(1)}(t;s,x)\in B_{\kappa,U}(t,x), \qquad s\le t. \label{eq:x1} \EEQ
 Then
\BEQ  \frac{d}{ds} x^{(2)}(t;s,x)=u_0(x^{(1)}\big(t-s,x^{(2)}(t;s,x))\big)\in 
u_0\big(B_{\kappa,U}(t,x^{(2)}(t;s,x)) \big).\EEQ

\medskip \noindent
This suggests considering generalizations of the flow $t\mapsto \Phi_{\kappa,U}(t,x)$
of the following kind:

\begin{Definition}[generalized flow] \label{def:generalized-flow}
Let $t,x_{min}>0$ and $\kappa>1$, $\tilde{U}\ge 1$. A {\em generalized flow} with {\em initial velocity}
$u_0$ and {\em parameters} $(\kappa,\tilde{U},x_{min})$ (in short, a $(\kappa,\tilde{U},x_{min})$-flow
with velocity $u_0$, or simply a $(\kappa,\tilde{U},x_{min})$-flow if $u_0$ is clear from
the context) is a system of ode's started from $x\in\R^d$,
\BEQ \frac{d}{ds}x(t;s,x)=u_0\big({\cal X}(t;s,x(t;s,x))\big), \qquad x(t;0,x)=x 
\label{eq:generalized-flow} \EEQ
with velocity field $v(t;s,\cdot)=u_0({\cal X}(t;s,\cdot))$ depending on the time-parameter $t$, such that
${\cal X}(t;s,y)\in B_{\kappa,\tilde{U},x_{min}}(t,y)$, $y\in\R^d$.

\noindent The mapping $(s,y)\mapsto {\cal X}(t;s,y)$ is simply called the
{\em mapping associated to the generalized flow} (\ref{eq:generalized-flow}).
\end{Definition}

Since our estimates concerning $(\kappa,\tilde{U},x_{min})$-flows do not depend
on $x_{min}$ provided $x_{min}\le (\tilde{U}t)^{\kappa/(\kappa-1)}$ (see Lemma
\ref{lem:BkappaU}), it is reasonable
to assume that $x_{min}\ge (\tilde{U}t)^{\kappa/(\kappa-1)}$ in the above Definition.

\medskip

{\bf In the sequel, $U$ is a fixed parameter associated to the growth at infinity of the
initial velocity $u_0$, while we let $\tilde{U}$ vary in some range included in
$[U,+\infty)$. }

\medskip

Under (Hyp1) such flows may be bounded very easily:


\begin{Lemma} \label{lem:generalized-flow-bound}
There exists some constant $C\ge 1$ such that the following holds. 
Let $\tilde{U}\ge U\ge 1$, $t\ge \tilde{U}^{-1}$, and ${\cal X}(t;.,.)$ be the mapping associated to a $(\kappa,\tilde{U},x_{min})$-flow. Assume
the initial velocity $u_0$ satisfies (Hyp1). Then $x(t;s,x)\in B_{\kappa,CU,C h_{\kappa}(t;\tilde{U},x_{min})}(t,x)$ for all $s\le t$, where
\BEQ 
h_{\kappa}(t,\tilde{U};x_{min}):=\tilde{U} t \left( \max(x_{min},(\tilde{U}t)^{\kappa/(\kappa-1)})
\right)^{1/\kappa}. \label{eq:g-kappa} \EEQ
\end{Lemma}

Of course, this result holds for arbitrary small $t$  provided one replaces $\tilde{U}t$
by $\langle \tilde{U}t\rangle$.
Note the particular case, 
\BEQ h_{\kappa}(t,\tilde{U};(\tilde{U}t)^{\kappa/(\kappa-1)})=
(\tilde{U}t)^{\kappa/(\kappa-1)}. \label{eq:gh} \EEQ

\medskip
\noindent {\bf Proof.} Clearly we may replace $x_{min}$ by $\max(x_{min},(\tilde{U}t)^{\kappa/(\kappa-1)})$. Hence we assume $x_{min}\ge (\tilde{U}t)^{\kappa/(\kappa-1)}$
is a large cut-off, 
and use Lemma \ref{lem:BkappaU} (2) in the following form,
\BEQ |{\cal X}(t;s,y)-y|\le C\tilde{U}t \max(x_{min},|y|)^{1/\kappa}. \EEQ

We distinguish two cases:
\begin{itemize}
\item[(i)] $(|y|\le x_{min})$  By (Hyp1) 
\BEA  |u_0({\cal X}(t;s,y))| &\le &  U\left(1+|y|+C\tilde{U}t\  x_{min}^{1/\kappa}\right)^{1/\kappa} \nonumber\\
&= &  U \left(1+ |y|+Ch_{\kappa}(t,U;x_{min}) \right)^{1/\kappa}; \label{eq:ok} \EEA

\item[(ii)] ($|y|\ge x_{min}$) By (Hyp1) again
\BEA  |u_0({\cal X}(t;s,y))| & \le &   U\left(1+|y|+C\tilde{U}t\  |y|^{1/\kappa}\right)^{1/\kappa}
\nonumber\\
&\le &    U C^{1/\kappa} (1+|y|+\tilde{U}t|y|^{1/\kappa})^{1/\kappa}
\le U(2C)^{1/\kappa} (1+|y|)^{1/\kappa}\le CU (1+|y|)^{1/\kappa}
\nonumber\\  \label{eq:better-than-ok} \EEA
for $C$ large enough;
\end{itemize}

which proves the Lemma.
    \hfill\eop

In particular we have proved: $x^{(2)}(t;s,x)\in B_{\kappa,CU,C(Ut)^{\kappa/(\kappa-1)}}(t,x)$ for all $s\le t$.  
\medskip

\noindent We may now iterate, and get for $m\ge 0$ and $t\ge U^{-1}$, using
(\ref{eq:gh}),
\BEQ x^{(m)}(t;s,x)\in B_{\kappa,CU, x_{min}^{(m)} }(t,x), \qquad s\le t\EEQ
with $x_{min}^{(0)}=x_{min}^{(1)}=0$, $x_{min}^{(2)}=C(Ut)^{\kappa/(\kappa-1)}$, and 
\BEQ x_{min}^{(m+1)}=Ch_{\kappa}(t,CU;x_{min}^{(m)})= C^2 Ut (x_{min}^{(m)})^{1/\kappa}, \qquad m\ge 2 \label{eq:Um}. \EEQ 
 This increasing recursive sequence converges for $m\to\infty$ for all $\kappa>1$; 
we get by Lemma  \ref{App:lem:alpha} a uniform bound for all $m\ge 0$,
\BEQ  x_{min}^{(m)}\le x_{min}^{(\infty)} \label{eq:Uinfty} \EEQ
where $x^{(\infty)}\lesssim (Ut)^{\kappa/(\kappa-1)}$
is the fixed point of the sequence.

\medskip
All this strongly suggests that the approximation scheme should converge under the
hypothesis (Hyp1).  
Leaving any rigor at this stage, and letting $m\to\infty$, one may conjecture
that the solution of Burgers' equation satisfies for $t\ge U^{-1}$
\BEQ u(t,x)\in u_0\left(B_{\kappa,CU,C(Ut)^{\kappa/(\kappa-1)}}(t,x)\right). \label{eq:u-u0} \EEQ
Assuming (Hyp1), we get, using (\ref{eq:ok}) and (\ref{eq:better-than-ok}),
\BEQ |u(t,x)|\lesssim U (|x|+\langle Ut\rangle^{\kappa/(\kappa-1)})^{1/\kappa}. \label{eq:u-Hyp1} \EEQ
Note however that, contrary to (\ref{eq:Uinfty}), this bound strongly relies on (Hyp1).
When we consider later on more general initial conditions, (\ref{eq:u-Hyp1}) will be
replaced by a much weaker bound, see (\ref{eq:bound-u}) in Section 3.

\subsection{The viscous case}


We now come back to non-zero viscosity; we fix for simplicity $\eta=1$. 
Instead of (\ref{eq:xmxm-1}), we consider the approximation scheme (\ref{eq:um}) and
its Feynman-Kac solution (\ref{eq:FK1},\ref{eq:FK2}). To avoid dealing with stochastic
calculus tools we replace the stochastic differential equation (\ref{eq:FK2}) with
an ode with random coefficients by letting $Y^{(m)}(t;s,x):=X^{(m)}(t;s,x)-B_s$, a
conventional trick which is sometimes called the Doss-Sussmann trick: 
we thus get
\BEA \frac{d}{ds} Y^{(m)}(t;s,x) &=& u^{(m-1)}(t-s,Y^{(m)}(t;s,x)+B_s) \nonumber\\
&=& \tilde{\esper}\left[ u_0(X^{(m-1)}(t-s,Y^{(m)}(t;s,x)+B_s))\right] \nonumber\\
&=& \tilde{\esper}\left[ u_0(\tilde{B}_{t-s}+Y^{(m-1)}(t-s,Y^{(m)}(t;s,x)+B_s))\right] 
\label{eq:XmXm-1} \EEA
where $X^{(m-1)}(t-s,y)=\tilde{B}_{t-s}+Y^{(m-1)}(t-s,y)$ is a random characteristic
depending on an extra Wiener process $(\tilde{B}_t)_{t\ge 0}$, independent from $B$,
and $\tilde{\esper}[\, \cdot\, ]$ is the partial expectation with respect to $\tilde{B}$.
From standard results on Brownian motion,  $\sup_{0\le s\le t}|\tilde{B}_s|$ scales like
$\sqrt{t}$ and is actually bounded by $O(\sqrt{t})$ with high probability, namely,
there exists a constant $c>0$ such that $\proba[\sup_{0\le s\le t}|\tilde{B}_s|>A\sqrt{t}]\lesssim e^{-cA^2}$ for all $A>0$. In the ensuing discussion we introduce the rescaled
random variables,
\BEQ M_t:=1+ \frac{\sup_{0\le s\le t}|{B}_s|}{\sqrt{t}}, \qquad
\tilde{M}_t:=1+ \frac{\sup_{0\le s\le t}|\tilde{B}_s|}{\sqrt{t}}
 \label{eq:Mt} \EEQ
which are therefore $O(1)$ with high probability. In particular, for all 
$\alpha,A\ge 1$,
\BEQ \esper[(M_t)^{\alpha}]=O(1) \EEQ
with a constant depending on $\alpha$, 
\BEQ \proba[M_t>A]\lesssim e^{-cA^2} \EEQ
for some universal constant $c$, and similarly for $\tilde{M}_t$.

\medskip\noindent

\medskip
Let us consider for the sake of illustration the cases $m=0,1$. First

\BEQ Y^{(0)}(t;s,x)=x; \EEQ
solving explicitly the trivial $0$-th transport equation $(\partial_t-\Del)u^{(0)}(t,x)=0$, we get
\BEA \frac{d}{ds} Y^{(1)}(t;s,x) &=&\tilde{\esper}[u_0(Y^{(1)}(t;s,x)+B_s+\tilde{B}_{t-s})] \nonumber\\
&=& u^{(0)}(t-s;Y^{(1)}(t;s,x)+B_s)=e^{(t-s)\Del}u_0(Y^{(1)}(t;s,x)+B_s).
\label{eq:2.39} \EEA
It is easy to check that 
\BEQ e^{t\Del}(y\mapsto (1+|y|)^{1/\kappa})(x)\lesssim 1+t^{1/2\kappa}+|x|^{1/\kappa}.
\label{eq:etDely1/kappa} \EEQ
 Thus
\BEQ \big| \frac{d}{ds}Y^{(1)}(t;s,x)\big|\lesssim U\left(1+t^{1/2\kappa}+|B_s|^{1/\kappa}+|Y^{(1)}(t;s,x)|^{1/\kappa} \right). \label{eq:viscous1} \EEQ
 Note that the same result may be retrieved without solving for $u^{(0)}$: namely,
\BEA && \left| \tilde{\esper}[u_0(Y^{(1)}(t;s,x)+B_s+\tilde{B}_{t-s})] \right| \le
U \left\{ \tilde{\esper}\left[ 1+|\tilde{B}_{t-s}|+|B_s|+|Y^{(1)}(t;s,x)| \right] \right\}^{1/\kappa}  \nonumber\\
&& \qquad\qquad\qquad \le U\left(1+M_t\sqrt{t}+|Y^{(1)}(t;s,x)|\right)^{1/\kappa} \EEA
where we have used Jensen's inequality.

Hence (by definition)  $Y^{(1)}(t;s,x)\in B_{\kappa,U,1+M_t\sqrt{t})}$, implying in particular  
\BEQ Y^{(1)}(t;s,x)\in B_{\kappa,CU,\max(\langle Ut\rangle ^{\kappa/(\kappa-1)},M_t\sqrt{t})}(t,x),
\label{eq:Y1} \EEQ
with the advantage that the  cut-off is always {\em large} in this expression, in the
sense of Lemma \ref{lem:generalized-flow-bound} (2). We may distinguish two regimes:

\begin{itemize}
\item[(i)] $M_t\sqrt{t}>\langle Ut\rangle^{\kappa/(\kappa-1)}$ ({\em diffusion prevails over
convection}) then $Y^{(1)}(t;s,x)\in B_{\kappa,CU,M_t\sqrt{t}}(t,x)$, hence $|Y^{(1)}(t;s,x)-x|\lesssim Ut\ \max(|x|,M_t\sqrt{t})^{1/\kappa}$. 
\medskip

This case (i) is highly improbable if $U\gg 1$ (i.e. when convection effects
are important) since 
\BEQ \left(M_t\sqrt{t}\gtrsim \langle Ut\rangle ^{\kappa/(\kappa-1)} \right) \Longrightarrow
\left(M_t\gtrsim U^{1/2} \langle Ut\rangle^{\half \frac{\kappa+1}{\kappa-1}}\ge U^{1/2} \right)
\label{eq:U1/2} \EEQ
both if $t\le U^{-1}$ and $t\ge U^{-1}$.  
For $t$ large enough (depending on the random variable $M_t$) one is necessarily in case (ii);  

\item[(ii)]  $M_t\sqrt{t}\lesssim \langle Ut\rangle^{\kappa/(\kappa-1)}$ ({\em convection prevails
over diffusion}), then we simply get
$Y^{(1)}(t;s,x)\in B_{\kappa,C'U,(C'\langle Ut\rangle) ^{\kappa/(\kappa-1)}}(t,x)$.
\end{itemize}

\medskip \noindent As in the non-viscous case, we want to iterate. To go further, we
need a rather straightforward adapatation to the viscous case of the notion of 
generalized $(\kappa,\tilde{U},x_{min})$-flow introduced in the previous paragraph.

\begin{Definition}[viscous generalized flow] 
\label{def:viscous-generalized-flow} (compare with Definition \ref{def:generalized-flow}) 
Let $>0$ and $\kappa>1$, $\tilde{U}\ge 1$. A {\em viscous generalized flow} with {\em initial velocity}
$u_0$ and {\em parameters} $(\kappa,\tilde{U},X_{min})$ (in short, a viscous $(\kappa,\tilde{U},X_{min})$-flow
with velocity $u_0$, or simply a viscous $(\kappa,\tilde{U},X_{min})$-flow if $u_0$ is clear from
the context) is a system of ode's with random coefficients started from $x\in\R^d$,
\BEQ \frac{d}{ds}Y(t;s,x)=\tilde{\esper}\Big[u_0\big(\tilde{B}_{t-s}+ {\cal Y}(t;s,Y(t;s,x)+B_s))\Big], \qquad Y(t;0,x)=x 
\label{eq:viscous-generalized-flow} \EEQ
with random velocity field $v(t;s,\cdot)=\tilde{\esper}\big[u_0(\tilde{B}_{t-s}+ {\cal Y}(t;s,\cdot+B_s))\big]$ depending on the time-parameter $t$, such that
${\cal Y}(t;s,y)\in B_{\kappa,\tilde{U},X_{min}}(t,y)$, $y\in\R^d$, where $X_{min}=X_{min}(t)$ 
is a random variable depending on $(\tilde{B}_s)_{s\in[0,t]}$.

\noindent The mapping $(s,y)\mapsto {\cal Y}(t;s,y)$ is  called the
{\em mapping associated to the viscous generalized flow} (\ref{eq:viscous-generalized-flow}).
\end{Definition}

In the above example, see (\ref{eq:Y1}), $X_{min}=C\max((Ut)^{\kappa/(\kappa-1)},
M_t\sqrt{t})$.

Lemma \ref{lem:generalized-flow-bound} generalizes under (Hyp1) to the viscous case in the
following way.

\begin{Lemma} \label{lem:viscous-generalized-flow-bound}
There exists some constant $C\ge 1$ such that the following holds. 
Let $t\ge \tilde{U}^{-1}$ and ${\cal Y}(t;.,.)$ be the mapping associated to a viscous
generalized $(\kappa,\tilde{U},\max(x_{min},\tilde{M}_t\sqrt{t}))$-flow, with $x_{min}\ge (\tilde{U}t)^{\kappa/(\kappa-1)}$
{\em deterministic}. Assume
the initial velocity $u_0$ satisfies (Hyp1). Then 
\BEQ {\cal Y}(t;s,x)\in B_{\kappa,CU,\max(Ch_{\kappa}(t,\tilde{U};x_{min}),M_t\sqrt{t})}(t,x) \EEQ  for all $s\le t$, where
$h_{\kappa}(t,\tilde{U};x_{min}):=\tilde{U}t\,  x_{min}^{1/\kappa},$ as in Lemma
\ref{lem:generalized-flow-bound}.
\end{Lemma}

As in Lemma \ref{lem:generalized-flow-bound}, we note that this result holds for arbitrary small $t$ provided one replaces $\tilde{U}t$ by $\langle\tilde{U}t\rangle$.

Comparing with Lemma \ref{lem:generalized-flow-bound}, one sees that
the cut-off is larger due to diffusion in the
highly improbable regime, defined by $M_t\sqrt{t}>x_{min}$, where diffusion prevails over convection.
\medskip

{\bf Proof.} We distinguish two regimes:

\begin{itemize}
\item[(i)] $(|y+B_s|\le \max(x_{min},\tilde{M}_t\sqrt{t}))$. Then

\BEQ |u_0(\tilde{B}_{t-s}+{\cal Y}(t;s,y+B_s))|\le  U\left(1+\tilde{M}_t\sqrt{t}+ 
|y+B_s|+C  
\tilde{U}t \left( \max(x_{min}, \tilde{M}_t\sqrt{t}) \right)^{1/\kappa} \right)^{1/\kappa},\EEQ
whence (using $x_{min}/\sqrt{t}\ge (Ut)^{\kappa/(\kappa-1)}/\sqrt{t}\ge U^{1/2}\ge 1$,
see (\ref{eq:sqrt-t}))
\BEA && \Big| \tilde{\esper}\left[ {\bf 1}_{|y+B_s|\le \max(x_{min},\tilde{M}_t\sqrt{t})} u_0\big(\tilde{B}_{t-s}+ {\cal Y}(t;s,y+B_s)) \right] \Big|
\nonumber\\
&& \qquad \lesssim U
 \left(1+\sqrt{t}+ |y+B_s|+  \tilde{U}t (\max(x_{min},\sqrt{t}))^{1/\kappa}
 \right)^{1/\kappa} \nonumber\\
&& \lesssim U \left( |y|+M_t\sqrt{t}+
\tilde{U}t \ x_{min}^{1/\kappa} \right)^{1/\kappa} \label{eq:ok2} \EEA
as expected; 

\item[(ii)]  $(|y+B_s|\ge \max(x_{min},\tilde{M}_t\sqrt{t}))$. Then
\BEA  |u_0(\tilde{B}_{t-s}+{\cal Y}(t;s,y+B_s))| &\lesssim &   U\left(1+\tilde{M}_t\sqrt{t}+ 
|y+B_s|+  
\tilde{U}t |y+B_s|^{1/\kappa} \right)^{1/\kappa} \nonumber\\
&\lesssim & U|y+B_s|^{1/\kappa} \lesssim U (|y|+M_t\sqrt{t})^{1/\kappa}; \label{eq:better-than-ok2} \EEA 
which proves the Lemma.

\end{itemize}

 \hfill \eop

\medskip
Iterating as in the non-viscous case, we get for $m\ge 0$ and $t\ge U^{-1}$
\BEQ Y^{(m)}(t,x)\in B_{\kappa,CU,\max(x_{min}^{(m)}, M_t\sqrt{t})}(t,x) \EEQ
with as in the non-viscous case, see (\ref{eq:Um}), 
\BEQ x_{min}^{(0)}=x_{min}^{(1)}=0, \qquad x_{min}^{(2)}=C(Ut)^{\kappa/(\kappa-1)}, \qquad
x_{min}^{(m+1)}=Ch_{\kappa}(t,CU;x_{min}^{(m)})= C^2 Ut  (x_{min}^{(m)})^{1/\kappa} \label{eq:viscous-Um} \EEQ
bounded uniformly in $m$ by $O((Ut)^{\kappa/(\kappa-1)})$.

Assuming as in the non-viscous case that the approximation scheme converges, it is natural
to conjecture that the solution of Burgers' equation satisfies, still under (Hyp1)
\BEA |u(t,x)| &=& \Big|\lim_{m\to\infty} \esper[u_0(X^{(m)}(t,x))]  \Big| \nonumber\\
&\lesssim & U\  \esper\left[ \left(|x|+M_t\sqrt{t}+ (Ut)^{\kappa/(\kappa-1)}\right)^{1/\kappa} \right] \nonumber \\
&\lesssim & U  (|x|+(Ut)^{\kappa/(\kappa-1)})^{1/\kappa}   \nonumber\\  \label{eq:bound-u-Hyp1}
\EEA
(see proof of Lemma \ref{lem:viscous-generalized-flow-bound})
as in the non-viscous case.
\bigskip


\section{More general initial data}


From the previous section, in particular, Lemmas \ref{lem:generalized-flow-bound} and
\ref{lem:viscous-generalized-flow-bound}, it is reasonable to expect that the sequence
$(u^{(m)})_{m\ge 0}$ is controlled as soon as flows driven by $u_0$, or the 'generalized
flows' thereof introduced in Definition \ref{def:generalized-flow}, \ref{def:viscous-generalized-flow}, are controlled well enough, in particular for $t$ large,  so as to ensure the possibility of
an induction. This opens the way to flows subject to sudden but brief accelerations,
corresponding to small areas where $u_0$ may be indeed very large; those must be
brief enough so as not to change the behaviour of the flow for $t$ large. 
What 'large' means is not so clear. Here we are interested in the whole regime
$t\in[\frac{1}{U},+\infty)$.  

It would be  natural to think of {\em defining}
 $u_0$ to be {\em admissible} if Lemmas \ref{lem:generalized-flow-bound} and
\ref{lem:viscous-generalized-flow-bound}, or some generalization thereof, hold. We did not find
however any class of examples of admissible initial velocities $u_0$ which do not satisfy
(Hyp1). Instead, we shall construct  in the following way explicit examples of initial velocities for which we get uniform a priori
bounds for the characteristics. First we  consider some  $\tilde{u}_0$ satisfying
(Hyp1). Then we modify it in an essentially arbitrary way in a region with small relative volume,
from which it can therefore escape in arbitrarily short time. The main challenge is to prove that
there exist {\em safe zones}, with relative volume tending to 1 at spatial infinity, which
are essentially {\em stable} under the flows -- deterministically in the non-viscous case,
with high probability in the viscous case. This {\em safe zone stability property}
(see Theorem \ref{th:admissible}, Theorem \ref{th:viscous-admissible}) must be
proved by induction. Then the complementary of the safe zones is made of small, widely
separated islands, called {\em dangerous zones}, which by the safe zone stability property cannot communicate with each other;
this simple fact settles non-inductively the analysis of trajectories started outside safe zones.

Let us mention that  for a given velocity $u_0$ such that the associated flow has a 
relatively simple  large scale topological structure (including large limit cycles, etc.)
is not too complicated, the existence of large safe zones should not be too complicated to verify if true. Thus criteria (\ref{hyp:R1},\ref{hyp:R2}) below should merely be considered
as some option.

\begin{Definition} \label{def:dangerous-zones}
Let $(R_n)_{n\ge 1}$ be an increasing sequence, $1\le R_1<R_2<R_3<\ldots$ such that, for all $i\ge 1$, 
\BEQ R_{2i}-R_{2i-1}\le R_{2i-1}^{1/\kappa}, \label{hyp:R1} \EEQ
\BEQ  R_{2i+1}\ge 4 R_{2i} \label{hyp:R2}
.\EEQ
Annuli $B(0,R_{2i+1})\setminus B(0,R_{2i})$ are called {\em safe zones}. Annuli
${\cal A}_i:=B(0,R_{2i})\setminus B(0,R_{2i-1})$ are called {\em dangerous zones}. 
\end{Definition}

\begin{Remark} \label{rem:subdivide} For  convenience we repeatedly subdivide any large safe zone $B(0,R_{2i+1})\setminus B(0,R_{2i})$ such that
$R_{2i+1}\ge 16R_{2i}$ into $\Big(B(0,4R_{2i})\setminus B(0,R_{2i}) \Big) \uplus \emptyset
\uplus \Big( B(0,R_{2i+1})\setminus B(0,4R_{2i})\Big)$, with an empty dangerous zone
sandwiched in-between, until all safe zones $B(0,R_{2i+1})\setminus B(0,R_{2i})$ are such that
$R_{2i+1}< 16 R_{2i}$. 
\end{Remark}

As explained in the introduction, our results hold if $R_{2i}-R_{2i-1}\le CR_{2i-1}^{1/\kappa}$ and $R_{2i+1}\ge (1+\eps)R_{2i}$ for some $C,\eps>0$. We imposed (\ref{hyp:R1},\ref{hyp:R2}) because we did not want to make explicit the dependence of our
bounds on $C,\eps$. 

\medskip

\noindent We first consider the simpler non-viscous case.


\subsection{Non-viscous case}


To give a flavor of the proofs of Theorems \ref{th:admissible} and \ref{th:viscous-admissible}
below, we start with the following elementary Lemma. It helps choosing a constant $C>1$ such
that  \BEQ |y-x|\le (C-1)Ut \max(x_{min},|x|)^{1/\kappa} \label{eq:C-1} \EEQ
provided  $x_{min}\ge (Ut)^{\kappa/(\kappa-1)}$ and $y\in B_{\kappa,U,x_{min}}(t,x)$ (see Lemma \ref{lem:BkappaU} (2)). 

In order to take into account various numerical constants coming from elementary estimates 
(Taylor expansions, etc.), {\bf we assume once and for all that
$C$ is large enough}. 

\medskip

\begin{Lemma} \label{lem:admissible}
Let $\tilde{u}_0:\R^d\to\R^d$ be an initial $C^1$ velocity satisfying (Hyp1) for some constants
$U\ge 1$, $\kappa>1$. Let $u_0:\R^d\to\R^d$ be any Lipschitz function coinciding with
$\tilde{u}_0$ outside the union of annuli $\cup_{i\ge 1} {\cal A}_i$, ${\cal A}_i:=B(0,R_{2i})\setminus B(0,R_{2i-1})$. Then the solution of the ode 
$\frac{dy}{ds}=u_0(y)$, $y(0)=x$, satisfies
\BEQ |y(s)-x|\le 16 (C-1) \langle Ut\rangle  \max((16C\langle Ut\rangle)^{\kappa/(\kappa-1)},|x|)^{1/\kappa}, \qquad 0\le s\le t. \label{eq:admissible} \EEQ
\end{Lemma}

{\bf Proof.} 

Let us first make a general remark. If $u_0\equiv\tilde{u}_0$ along the whole trajectory $(y(s))_{0\le s\le t}$, then
$y(s)$ is bounded as in (\ref{eq:C-1}), where we have set $x_{min}=(Ut)^{\kappa/(\kappa-1)}$,
\BEQ |y(s)-x|\le (C-1) Us \max((Ut)^{\kappa/(\kappa-1)},|x|)^{1/\kappa}.\EEQ

We must now distinguish two cases.

\begin{itemize}
\item[(i)] Let $|x|\ge (16 C\langle Ut\rangle)^{\kappa/(\kappa-1)}$ (later on we shall actually need to assume
that $|x|\ge 32(16C \langle Ut\rangle)^{\kappa/(\kappa-1)}$).  Then $|x|^{1/\kappa}\le |x|/16CUt$, so,
{\em provided} $u_0\equiv \tilde{u}_0$ {\em along the whole trajectory},
\BEQ |y(s)|\ge |x|-(C-1)Ut|x|^{1/\kappa} \ge \frac{|x|}{2}, \qquad 
|y(s)|\le |x|+ (C-1)Ut|x|^{1/\kappa} \le 2|x|.  \label{eq:along} \EEQ
Thus we check a posteriori that $u_0\equiv \tilde{u}_0$ along the whole trajectory if
\BEQ |x|\in I_i(t):=[R_{2i}+4(C-1)Ut R_{2i}^{1/\kappa}, R_{2i+1}-4(C-1)UtR_{2i+1}^{1/\kappa}]  \label{eq:|x|} \EEQ
(with $C$ large enough as stipulated above), with $R_{2i}\ge (16 C\langle Ut\rangle)^{\kappa/(\kappa-1)}$; note that {\em if} $|x|\ge 16 (16CUt)^{\kappa/(\kappa-1)}$ and  (\ref{eq:|x|}) holds, then indeed $R_{2i}\ge \frac{1}{16} R_{2i+1}\ge (16C\langle Ut
\rangle)^{\kappa/(\kappa-1)}$ by construction. Namely, if $|x|\in I_i(t)$ then 
\BEA && \label{eq:Imin}  (R_{2i+1}-4(C-1)UtR_{2i+1}^{1/\kappa}) + (C-1)Us \left(R_{2i+1}-4(C-1)UtR_{2i+1}^{1/\kappa}
\right)^{1/\kappa} \nonumber\\
&& \qquad \qquad \le R_{2i+1}-4(C-1)U(t-s)R_{2i+1}^{1/\kappa}; \EEA

\BEA
&& (R_{2i}+4(C-1)UtR_{2i}^{1/\kappa}) - (C-1)Us \left(R_{2i}+4(C-1)UtR_{2i}^{1/\kappa}
\right)^{1/\kappa}  \nonumber\\
&&\qquad \qquad \ge R_{2i}+4(C-1)UtR_{2i}^{1/\kappa}-(C-1)Us (2R_{2i})^{1/\kappa}  \nonumber\\
&& \qquad \qquad\ge R_{2i}+4(C-1)U(t-s)R_{2i}^{1/\kappa}  
\label{eq:Imax}
\EEA
so 
\BEQ |y(s)|\in I_i(t-s)\subset I_i(0)=[R_{2i},R_{2i+1}]. \label{eq:stable} \EEQ

We call $(I_i(t))_i$  {\em safe intervals}; (\ref{eq:stable}) is the main argument in our
{\em safe zone stability property}.  Note that $I_i(t)\not=\emptyset$ since
\BEQ (R_{2i+1}-4(C-1)Ut R_{2i+1}^{1/\kappa})-(R_{2i}+4(C-1)UtR_{2i}^{1/\kappa})\ge 
\half R_{2i+1}-\frac{3}{2} R_{2i} \ge \half R_{2i} \EEQ
by Hypothesis (\ref{hyp:R2}).

\medskip
If now $x$ does not belong to a safe zone, say, $|x|\in [R_{2i-1}-4(C-1)Ut R_{2i-1}^{1/\kappa},
R_{2i}+4(C-1)UtR_{2i}^{1/\kappa}]$, then $x$ is possibly free to move in essentially arbitrarily small time to
$x'=y(t')$, $t'\in[0,t]$, such that $|x'|$ is the closest end of one of the two neighbouring
safe zones, $I_j(t)$, with $j=i-1$ or $i$.  Then for $C$ large enough we get successively, using as unique ingredients Hypotheses (\ref{hyp:R1},\ref{hyp:R2}) and the lower bound 
$|x|\ge 16 (16 C \langle Ut\rangle)^{\kappa/(\kappa-1)}$, 
$$R_{2i}\ge (8C\langle Ut\rangle)^{\kappa/(\kappa-1)};  $$
$$ R_{2i-1}\ge R_{2i}-O(R_{2i}^{1/\kappa})\ge \frac{3}{4} R_{2i}\ge (6C\langle Ut\rangle)^{\kappa/(\kappa-1)}; $$
$$|x|\ge R_{2i-1}-4(C-1)UtR_{2i-1}^{1/\kappa} \ge \frac{R_{2i-1}}{3} \ge \frac{R_{2i}}{4}; $$
\BEA |x'-x| &\le & (R_{2i}+4(C-1)Ut R_{2i}^{1/\kappa})-(R_{2i-1}-4(C-1)Ut R_{2i-1}^{1/\kappa}) 
\nonumber\\
&\le & (C-1)  R_{2i-1}^{1/\kappa}  \left\{1+2(C-1) Ut\,  R_{2i-1}^{(1/\kappa)-1} 
\right\} \nonumber\\
&\le & 2(C-1) R_{2i-1}^{1/\kappa}; \EEA
\BEA |x'| & \ge & R_{2i}-|x'-x|\ge  R_{2i}- 2CR_{2i-1}^{1/\kappa} \ge R_{2i}-2CR_{2i}^{1/\kappa} \nonumber\\
& \ge & \frac{3}{4} R_{2i}\ge \frac{1}{2}(R_{2i}+4(C-1)UtR_{2i}^{1/\kappa}) \ge \frac{|x|}{2}. \label{eq:x'>} \EEA

Assume $|x|\ge 32 (16C\langle Ut\rangle)^{\kappa/(\kappa-1)}$. If $j=i$ then $|x'|\ge |x|$ and
$I_j(t)\subset[(16C\langle Ut\rangle)^{\kappa/(\kappa-1)},\infty)$; otherwise 
$\min(I_j(0))\ge \frac{1}{16} \max(I_j(t))=\frac{|x'|}{16}\ge \frac{|x|}{32}$, so we get the
same conclusion. Thus the rest of the trajectory (for $s\ge t'$) remains inside a safe zone and
(\ref{eq:along}) holds, $|y(s)-x'|\le (C-1)Ut\, |x'|^{1/\kappa}$.  Hence for every $s\in[0,t]$, we get
\BEA  |y(s)-x| &\le&  2(C-1)R_{2i-1}^{1/\kappa}+ (C-1)Ut (R_{2i}+4(C-1)UtR_{2i}^{1/\kappa})^{1/\kappa}  \nonumber\\
&\le&  4(C-1) Ut R_{2i}^{1/\kappa}
\le 16 (C-1) Ut |x|^{1/\kappa}. \label{eq:3.11} \EEA
Note that (\ref{eq:3.11}) improves on (\ref{eq:admissible}) in the initial time regime $Ut\le 1$.

\item[(ii)] Let $|x|\le 32(16C\langle Ut\rangle)^{\kappa/(\kappa-1)}$. Then either the whole trajectory is
contained in $B(0,32(16C\langle Ut\rangle)^{\kappa/(\kappa-1)}$, or, letting $t'=\inf\{s\in[0,t]\ |\ |y(s)|=32(16C\langle Ut\rangle)^{\kappa/(\kappa-1)}\}$, we get by (i)
\BEQ |y(s)-y(t')|\le 16 (C-1) Ut |y(t')|^{1/\kappa} \le 32 (16C\langle Ut\rangle)^{\kappa/(\kappa-1)}, \qquad 
s\in[t',t] \EEQ
hence in whole generality, $|y(s)-x|\le 96(16C\langle Ut\rangle)^{\kappa/(\kappa-1)}$, $s\in [0,t]$. 

\end{itemize} \hfill \eop

\medskip
Now comes the main result.

\begin{Theorem}[non-viscous case] \label{th:admissible}
Let $\tilde{u}_0:\R^d\to\R^d$ be an initial $C^1$ velocity satisfying (Hyp1) for some constants
$U\ge 1$, $\kappa>1$. Let $u_0:\R^d\to\R^d$ be any Lipschitz function coinciding with
$\tilde{u}_0$ outside the union of annuli $\cup_{i\ge 1} {\cal A}_i$, ${\cal A}_i:=B(0,R_{2i})\setminus B(0,R_{2i-1})$. Then  the sequence of characteristics
$(x^{(m)}(t;\cdot,x))_{m\ge 0}$ satisfies the following uniform in $m$ estimates:

\begin{itemize}
\item[(i)] Let $|x|\ge (16C \langle Ut\rangle)^{\kappa/(\kappa-1)}$, then $|x^{(m)}(t;s,x)-x|\lesssim
(C-1)Ut|x|^{1/\kappa}$. If furthermore $x$ is in a safe zone, $|x|\in I_i(t)$, such that
$I_i(0)\subset [(16C\langle Ut\rangle)^{\kappa/(\kappa-1)},\infty)$, then $|x^{(m)}(t;s,x)|\in I_i(t-s)$ for $0\le s\le t$  ({\em safe
zone stability property}).
\item[(ii)] Let $|x|\le (16C\langle Ut\rangle)^{\kappa/(\kappa-1)}$. Then $|x^{(m)}(t;s,x)-x|\lesssim 
(16C\langle Ut\rangle)^{\kappa/(\kappa-1)}$.  
\end{itemize}
\end{Theorem}

Note that these estimates have just been proved in the case $m=1$.  We subdivide the {\bf proof} into three points.

\begin{itemize}
\item[(1)] The core of the proof is the safe zone stability property. Let $i\ge 1$ such that $I_i(0)\subset [(16C\langle Ut\rangle)^{\kappa/(\kappa-1)},\infty)$. Assume by induction that (see (\ref{eq:along},\ref{eq:stable}))
\BEQ (|x|\in I_i(t))\Longrightarrow \left( |x^{(m-1)}(t,x)|\in I_i(0), \ \frac{|x|}{2}\le 
|x^{(m-1)}(t,x)|\le 2|x| \right). \label{eq:ind-hyp} \EEQ
For such an $x$, we therefore know that in the ode for $x^{(m)}(t;\cdot,x)$,
\BEQ \frac{d}{ds}y(s)=u_0(x^{(m-1)}(t-s,y(s))), \EEQ
the norm of the argument of $u_0$, $x^{(m-1)}(t-s,y(s))$, belongs to $I_i(0)$ provided $|y(s)|\in I_i(t-s)$.
If this is the case, then
\BEQ |\frac{d}{ds}y(s)|=|\tilde{u}_0(x^{(m-1)}(t-s,y(s))|\le U(1+2|y(s)|)^{1/\kappa}\le 2^{1/\kappa}
U(1+|y(s)|)^{1/\kappa} \EEQ
by our induction hypothesis (\ref{eq:ind-hyp}), hence
\BEQ \frac{|x|}{2}\le |x|-2^{1/\kappa}(C-1)Ut|x|^{1/\kappa} \le |y(s)|\le |x|+2^{1/\kappa}
(C-1)Ut|x|^{1/\kappa}\le 2|x|. \EEQ
This leads to a slight modification of (\ref{eq:Imin},\ref{eq:Imax}),
\BEA && \label{eq:Imin-bis}  (R_{2i+1}-4(C-1)UtR_{2i+1}^{1/\kappa}) +2^{1/\kappa} (C-1)Us \left(R_{2i+1}-4(C-1)UtR_{2i+1}^{1/\kappa}
\right)^{1/\kappa} \nonumber\\
&& \qquad \qquad \le R_{2i+1}-4(C-1)U(t-s)R_{2i+1}^{1/\kappa}; \EEA

\BEA
&& (R_{2i}+4(C-1)UtR_{2i}^{1/\kappa}) - 2^{1/\kappa}(C-1)Us \left(R_{2i}+4(C-1)UtR_{2i}^{1/\kappa}
\right)^{1/\kappa}  \nonumber\\
&&\qquad \qquad \ge R_{2i}+4(C-1)UtR_{2i}^{1/\kappa}-2^{1/\kappa}(C-1)Us (2R_{2i})^{1/\kappa}  \nonumber\\
&& \qquad \qquad\ge R_{2i}+4(C-1)U(t-s)R_{2i}^{1/\kappa}  
\label{eq:Imax-bis}
\EEA

Hence we have checked a posteriori the safe zone stability property, $|y(s)|\in I_i(t-s)$.

\medskip
\item[(2)]
Assume now $|x|\ge 32(16C\langle Ut\rangle)^{\kappa/(\kappa-1)}$ but $x$ does not belong to a safe zone,
say, $|x|\in[R_{2i-1}-4(C-1)UtR_{2i}^{1/\kappa},R_{2i}+4(C-1)UtR_{2i}^{1/\kappa}]$. From the
proof of Lemma \ref{lem:admissible}, we know that the trajectory, if ever, enters a safe interval
$I_j(t)$, $j=i-1$ or $i$, at some point $x'=y(t')$ such that $|x'|\ge \frac{|x|}{2}$, and
$I_j(0)\subset[(16C\langle Ut\rangle)^{\kappa/(\kappa-1)},\infty)$. Hence we can avail ourselves of the
safe zone stability property proved in (1), yielding $|y(t)-x'|\le 2^{1/\kappa}(C-1)Ut |x|^{1/\kappa}$. Thus, for all $s\in[0,t]$,
\BEA  |y(s)-x| &\le&  2(C-1)R_{2i-1}^{1/\kappa}+ 2^{1/\kappa} (C-1)Ut (R_{2i}+4(C-1)UtR_{2i}^{1/\kappa})^{1/\kappa}  \nonumber\\
&\le&  5(C-1) Ut R_{2i}^{1/\kappa}
\le 20 (C-1) Ut |x|^{1/\kappa},\EEA
as in (\ref{eq:3.11}), up to a different numerical constant.

\item[(3)]
Finally, for $|x|\le 32 (16C\langle Ut\rangle)^{\kappa/(\kappa-1)}$, we conclude as in point (ii) of the proof 
of Lemma \ref{lem:admissible}, again up to different numerical constants.

\end{itemize}
 \hfill \eop

Under the hypotheses of Theorem \ref{th:admissible},  we obtain as in the previous section a
conjectural uniform bound for $u^{(m)}$ and for $u$, which we write down for $u$,
\BEQ u(t,x)\in u_0(B_{\kappa,CU,C\langle Ut\rangle^{\kappa/(\kappa-1)}}(t,x)) \EEQ
for some constant $C$, see (\ref{eq:u-u0}),  which is however not as explicit as (\ref{eq:u-Hyp1}).


\subsection{Viscous case}


Let us now consider the viscous case.
\medskip

\noindent The new difficulty here is that, for $M_t\sqrt{t}$ or $\tilde{M}_t\sqrt{t}$ large,
we clearly lose our safe zone stability property. Hence we need some general  a priori
bound on $u_0$; a polynomial bound at infinity is a very weak but sufficient requirement. 
Apart from that, the scheme follows closely that of \S 3.1.

\begin{Lemma} \label{lem:viscous-admissible}
Let $\tilde{u}_0:\R^d\to\R^d$ be an initial velocity satisfying (Hyp1) for some constants
$U\ge 1$, $\kappa>1$. Let $u_0:\R^d\to\R^d$ be any Lipschitz function coinciding with
$\tilde{u}_0$ outside the union of annuli $\cup_{i\ge 1} {\cal A}_i$, ${\cal A}_i:=B(0,R_{2i})\setminus B(0,R_{2i-1})$ {\em and satisfying the following a priori bound}
for some constants $\alpha,\beta\ge 0$,
\BEQ |u_0(x)|\le K_0(1+|x|)^{\frac{\alpha}{2}+\frac{1}{\kappa}}, \qquad
x\in\R^d  \label{eq:a-priori-bound-u0}
\EEQ
with
\BEQ K_0\le U^{\frac{\beta}{2}+1}.\EEQ  
Then the solution of the ode 
$\frac{d}{ds} Y(s)=\tilde{\esper}\left[u_0(Y(s)+B_s+\tilde{B}_{t-s})\right]$, $Y(0)=x$
 (see (\ref{eq:2.39})), satisfies
\BEQ |Y(s)-x|\lesssim  (C-1) \langle U t\rangle \max\left(
(16C\langle Ut\rangle)^{\kappa/(\kappa-1)},|x| \right)^{1/\kappa},\EEQ
if  $M_t\sqrt{t}\le \max\left( \langle Ut\rangle^{\kappa/(\kappa-1)}, \langle Ut \rangle \langle x\rangle^{1/\kappa} \right)$,
\BEQ |Y(s)-x|\lesssim \left( \frac{M_t\sqrt{t}}{\langle Ut\rangle} \right)^{\kappa} \EEQ
if
 $M_t\sqrt{t}\ge \max\left(  \langle Ut\rangle^{\kappa/(\kappa-1)},  \langle Ut\rangle  \langle x\rangle^{1/\kappa} \right)$.

\end{Lemma}

\medskip

 The {\bf proof} is a generalization of the non-viscous case, see proof of 
Lemma \ref{lem:admissible}.  We distinguish two regimes, (i) the {\em normal regime} where
{\em convection prevails over diffusion} ($M_t\sqrt{t}$ {\em small}), and (ii) the regime where {\em diffusion prevails
over convection} ($M_t\sqrt{t}$ large). The general idea is that the safe zone stability property holds in case (i),
while the a priori bound (\ref{eq:a-priori-bound-u0}) on $u_0$ yields new estimates in case (ii).
 Mind however (\ref{eq:a-priori-bound-u0}) is also needed in case (i) since $\tilde{M_t}\sqrt{t}$ may
be large. In particular (since a priori bounds alone would lead to a finite time explosion
of the paths), $|y|,U,t$  are controlled either {\em deterministically}
by $M_t$ -- which is not averaged over here -- or {\em stochastically} by $\tilde{M}_t$, when these get abnormally large.

\medskip

As usual, we may in practice assume that $|x|\ge (16C\langle Ut\rangle)^{\kappa/(\kappa-1)}$.

\begin{itemize}
\item[(i)] (normal regime) Assume $M_t\sqrt{t}\le \langle Ut\rangle \langle x\rangle^{1/\kappa}$. 
We first need an a priori bound of
\BEQ I^{(1)}:=\Big|\tilde{\esper}\left[ {\bf 1}_{\tilde{M}_t\sqrt{t}\ge \langle Ut\rangle \langle x\rangle^{1/\kappa}} u_0(Y(s)+B_s+
\tilde{B}_{t-s})\right] \Big|.\EEQ 
The event $\tilde{\Omega}:\tilde{M}_t\sqrt{t}\ge \langle Ut\rangle \langle x\rangle^{1/\kappa}$ is a rare even of probability
$O\left(\exp \ -c\left( \frac{\langle Ut\rangle}{\sqrt{t}}  \langle x\rangle^{1/\kappa} \right)^2 \right)=O(e^{-cUt}  e^{-cU\langle x\rangle^{2/\kappa}})$ (the last equality holds both for $Ut\le 1$ and $Ut\ge 1$!); thus $|x|$, but also $U$ and $t$, are 'stochastically' controlled
by $\tilde{M}_t$ (see below).
  {\em Provided} $|Y(s)|\lesssim |x|$ we get
 \BEQ I^{(1)}\lesssim K_0\  \tilde{\esper}\left[ {\bf 1}_{\tilde{\Omega}} (|x|+\tilde{M}_t\sqrt{t})^{\frac{\alpha}{2}+\frac{1}{\kappa}}\right].\EEQ
All factors in the above expression are highly suppressed by the exponentially small
factors $O(e^{-cUt} e^{-cU
\langle x\rangle^{2/\kappa}})$ since

\BEQ K_0\le  U^{\frac{\beta}{2}+1}\le (U\langle x\rangle ^{2/\kappa})^{\frac{\beta}{2}+1},
\qquad |x|\lesssim (U\langle x\rangle^{2/\kappa})^{\kappa/2}, \qquad \sqrt{t}\le (Ut)^{1/2}.\EEQ

Partitioning the event ${\bf 1}_{\tilde{M}_t\sqrt{t}\ge \langle Ut\rangle \langle x\rangle^{1/\kappa}} $
into $\cup_{n\ge 0} \tilde{\Omega}_n$ where $\tilde{\Omega}_n:=\{2^n \langle Ut\rangle 
\langle x\rangle^{1/\kappa} \le 
\tilde{M}_t\sqrt{t}< 2^{n+1} \langle Ut\rangle \langle x\rangle^{1/\kappa}  \}$, one can easily prove that,
 for $c'$ small enough,
\BEQ I\lesssim \sum_{n\ge 0} e^{-c'2^{2n} U\langle x\rangle^{2/\kappa}}
\lesssim 1.\EEQ

Hence, {\em provided} $|Y(s)|\approx |x|$, 
\BEQ \frac{d}{ds}Y(s)=O(1)+u_0(Y(s)+O(\langle Ut\rangle \langle x\rangle^{1/\kappa})=O(1)+u_0(Y(s)+O(\langle Ut\rangle \langle Y(s)\rangle^{1/\kappa}).\EEQ
Now, the innocuous replacement $Y(s)\mapsto Y(s)+O(\langle Ut\rangle \langle Y(s)\rangle^{1/\kappa})$ leaves the analysis of
Lemma \ref{lem:admissible} unchanged, up to the following modifications: define
\BEQ I_i(t):=[R_{2i}+2(C-1)(\langle Ut\rangle+Ut) R_{2i}^{1/\kappa}, R_{2i+1}-2(C-1)(\langle Ut\rangle+Ut)R_{2i+1}^{1/\kappa}] \EEQ
(compare with (\ref{eq:|x|})), so that the image of $I_i(0)=[R_{2i}+2(C-1)R_{2i}^{1/\kappa}, R_{2i+1}-2(C-1)R_{2i+1}^{1/\kappa}]$ by the mapping $y\mapsto y+O(\langle Ut\rangle \langle y\rangle^{1/\kappa})$ is $\subset[R_{2i},R_{2i+1}]$. For $C$ large enough
and $|x|\in I_i(t)$, one gets $Y(s)\in[R_{2i}+2(C-1)(\langle Ut\rangle+U(t-s)) R_{2i}^{1/\kappa}, R_{2i+1}-2(C-1)(\langle Ut\rangle+U(t-s))R_{2i+1}^{1/\kappa}]\subset I_i(t-s)$.

\item[(ii)] Assume on the contrary $M_t\sqrt{t}\ge \langle Ut \rangle \langle x\rangle^{1/\kappa}$; thus
$\langle x\rangle$ is controlled in terms of $M_t\sqrt{t}$,
\BEQ \langle x\rangle\le \left( \frac{M_t\sqrt{t}}{\langle Ut\rangle } \right)^{\kappa}\le (M_t\sqrt{t})^{\kappa}. \EEQ
 Thus the bound for $I^{(1)}$, see (i), is modified as follows 
{\em provided} $\langle Y(s)\rangle \le \left( \frac{M_t\sqrt{t}}{\langle Ut\rangle} \right)^{\kappa}\le (M_t\sqrt{t})^{\kappa}$,
\BEA I^{(1)} &\lesssim & K_0\  \tilde{\esper} \left[ {\bf 1}_{\tilde{\Omega}} (|Y(s)|+\tilde{M}_t\sqrt{t}
+M_t\sqrt{t})^{\frac{\alpha}{2}+\frac{1}{\kappa}} \right] \nonumber\\
&\lesssim & K_0 \left( |Y(s)|^{\frac{\alpha}{2}+\frac{1}{\kappa}} + 
\tilde{\esper}[ {\bf 1}_{\tilde{\Omega}} (\tilde{M}_t\sqrt{t})^{\frac{\alpha}{2}+\frac{1}{\kappa}} ]
+ (M_t\sqrt{t})^{\frac{\alpha}{2}+\frac{1}{\kappa}} \right) \nonumber\\
&\lesssim & U^{\frac{\beta}{2}+1} \max(1,(M_t\sqrt{t})^{\kappa})^{\frac{\alpha}{2}+\frac{1}{\kappa}}<\infty,
\label{eq:lem-I}
\EEA
to which one must add a smaller term,
\BEQ I^{(1),c}:=\Big|\tilde{\esper}\left[ {\bf 1}_{\tilde{\Omega}^c} u_0(Y(s)+B_s+
\tilde{B}_{t-s})\right] \Big| \lesssim K_0 (|Y(s)|^{\frac{\alpha}{2}+\frac{1}{\kappa}} +
(M_t\sqrt{t})^{\frac{\alpha}{2}+\frac{1}{\kappa}}). \label{eq:lem-Ic} \EEQ 
Clearly (considering only powers of $M_t$ for $t$ fixed), these are very poor estimates of
the velocity when $\frac{\alpha}{2}+\frac{1}{\kappa}>1$, given the a priori condition
$|Y(s)|=O(M_t^{\kappa})$; actually we shall not need them.

Now, it may happen that $\langle Y(t')\rangle=\left( \frac{M_t\sqrt{t}}{\langle Ut\rangle }\right)^{\kappa}$
 ($\ge |x|$) for some $t'\in[0,t]$. The estimates of (i) imply
 then in whole generality 
 \BEQ |Y(s)-x| \lesssim \left( \frac{M_t\sqrt{t}}{\langle Ut\rangle} \right)^{\kappa}.
 \label{eq:whole} \EEQ

\end{itemize}

 \hfill \eop

\medskip

We may now state the main theorem of this section, a counterpart of Theorem \ref{th:admissible}
in the viscous case. Safe intervals are defined as in the previous lemma.

\begin{Theorem}[viscous case] \label{th:viscous-admissible}
Let $\tilde{u}_0:\R^d\to\R^d$ be an initial velocity satisfying (Hyp1) for some constants
$U\ge 1$, $\kappa>1$. Let $u_0:\R^d\to\R^d$ be any Lipschitz function coinciding with
$\tilde{u}_0$ outside the union of annuli $\cup_{i\ge 1} {\cal A}_i$, ${\cal A}_i:=B(0,R_{2i})\setminus B(0,R_{2i-1})$, and satisfying the a priori bounds
(\ref{eq:a-priori-bound-u0}). Then  the sequence of characteristics
$(Y^{(m)}(t;\cdot,x))_{m\ge 0}$ satisfies the following uniform in $m$ estimates:

\begin{itemize}
\item[(i)] (normal regime, $M_t\sqrt{t}\le \max(\langle Ut\rangle ^{\kappa/(\kappa-1)},\langle Ut\rangle \langle x\rangle^{1/\kappa})$) 

Then $|Y^{(m)}(t;s,x)-x|\lesssim
(C-1)\langle Ut\rangle \max((16C\langle Ut\rangle)^{\kappa/(\kappa-1)}, |x|)^{1/\kappa}$. 

If furthermore $x$ is in a safe zone, $|x|\in I_i(t)$, such that
$I_i(t)\subset[(16C\langle Ut\rangle^{\kappa/(\kappa-1)},\infty)$, then, for all $x'\in\R^d$ such 
that $|x'-x|\le \langle Ut\rangle \langle x\rangle^{1/\kappa}$ and all $y\in\R^d$ such that $|y-Y^{(m)}(t;s,x')|\le 
\langle Ut\rangle \langle x\rangle^{1/\kappa}$, it holds $|y|\in I_i(t-s)$  ({\em safe
zone stability property}).

\item[(ii)] Assume $M_t\sqrt{t}\ge \max(\langle Ut\rangle^{\kappa/(\kappa-1)},
\langle Ut \rangle \langle x\rangle^{1/\kappa})$. Then
\BEQ |Y^{(m)}(t;s,x)-x|\lesssim \left( \frac{M_t\sqrt{t}}{\langle Ut\rangle}\right)^{\kappa}
.\EEQ
\end{itemize}

\end{Theorem}

{\bf Proof.} 
We proceed more of less as in the proof of Theorem \ref{th:admissible}. We cannot however  separate  the  inductive proof of the safe zone stability property from the rest of the argument since
we need the general bound (ii) to hold for $m-1$ to control the contribution to the velocity of
the event $\tilde{\Omega}:\tilde{M}_t\sqrt{t}\ge \max(Ut|x|^{1/\kappa},(Ut)^{\kappa/(\kappa-1)})$.
Thus we assume inductively that (i), (ii) hold for $m-1$. As usual, we may restrict the study to $|x|\ge (16C\langle Ut\rangle)^{\kappa/(\kappa-1)}$.

\begin{itemize}
\item[(i)] Assume first $M_t\sqrt{t}\le \max(\langle Ut\rangle ^{\kappa/(\kappa-1)}, \langle U\rangle \langle x\rangle^{1/\kappa})$
and let $x\in\R^d$ such that $|x|\in I_i(t)$, $I_i(t)\subset[(16C\langle Ut\rangle)^{\kappa/(\kappa-1)},\infty)$.
Recall $Y(s):=Y^{(m)}(t;s,x)$ solves the ode
\BEQ \frac{d}{ds} Y(s)=\tilde{\esper}\left[u_0(\tilde{B}_{t-s}+Y^{(m-1)}(t-s,Y(s)+B_s)) \right].\EEQ 
If $\tilde{M}_t\sqrt{t}\le \langle Ut\rangle \langle x\rangle^{1/\kappa}$, then $y:=\tilde{B}_{t-s}+Y^{(m-1)}(t-s,Y(s)+B_s)$
satisfies precisely the assumptions of the safe zone stability property, hence $|y|\in I_i(0)$
provided $|Y(s)|\in I_i(t-s)$.  Otherwise we first bound 
\BEQ I^{(m)}:=\left| \tilde{\esper}\left[ {\bf 1}_{\tilde{\Omega}} u_0(\tilde{B}_{t-s}+Y^{(m-1)}(t-s,Y(s)+B_s)) \right] \right|.\EEQ 
 Provided $|Y(s)|\approx |x|$ we get by induction hypothesis
\BEQ I^{(m)}\lesssim K_0\  \tilde{\esper}\left[ {\bf 1}_{\tilde{\Omega}} \left(|x|+\tilde{M}_t\sqrt{t}+
 \left( \frac{\tilde{M}_t\sqrt{t}}{\langle Ut\rangle}\right)^{\kappa}\right)^{\frac{\alpha}{2}+\frac{1}{\kappa}} \right] \lesssim 1 \EEQ
as in Lemma \ref{lem:viscous-admissible}. The  rest of the argument is as in the non-viscous case
(see proof of Theorem \ref{th:admissible}).

\item[(ii)] Assume now $M_t\sqrt{t}\ge \max(\langle Ut\rangle ^{\kappa/(\kappa-1)},
\langle Ut\rangle \langle x\rangle^{1/\kappa})$. 
By induction hypothesis we get
\BEA I^{(m)}  &\lesssim&  K_0\  \tilde{\esper}\left[ {\bf 1}_{\tilde{\Omega}}  \left(|Y(s)|
+M_t\sqrt{t}+\tilde{M}_t\sqrt{t}+
 \left( \frac{\tilde{M}_t\sqrt{t}}{\langle Ut\rangle}\right)^{\kappa}
 \right)^{\frac{\alpha}{2}+\frac{1}{\kappa}} \right] \nonumber\\
 &\lesssim & U^{\frac{\beta}{2}+1} \max(1,(M_t\sqrt{t})^{\kappa})^{\frac{\alpha}{2}+\frac{1}{\kappa}} <\infty \EEA
to which we must add a smaller contribution,
\BEA I^{(m),c} &:=& \left| \tilde{\esper}\left[ {\bf 1}_{\tilde{\Omega}^c} u_0(\tilde{B}_{t-s}+Y^{(m-1)}(t-s,Y(s)+B_s)) \right] \right| \nonumber\\
&\lesssim & K_0  \left(|Y(s)|
+M_t\sqrt{t}+
 \left( \frac{\tilde{M}_t\sqrt{t}}{\langle Ut\rangle}\right)^{\kappa}
 \right)^{\frac{\alpha}{2}+\frac{1}{\kappa}} 
 \EEA
 as in (\ref{eq:lem-I},\ref{eq:lem-Ic}). Using (i) one concludes as in (\ref{eq:whole}):
 $|Y(s)-x|\lesssim \left( \frac{M_t\sqrt{t}}{\langle Ut\rangle} \right)^{\kappa}.$
\end{itemize} \hfill \eop

\medskip

Using the  above Theorem we may conjecture that the following uniform bounds hold for
$u^{(m)}$, $m\ge 0$ and for $u$,

\BEA  |u(t,x)|&=& \left|\lim_{m\to\infty} \esper\left[u_0(X^{(m)}(t,x))\right] \right| \nonumber\\
&\lesssim & K_0 \esper\left[ \left( |x|+M_t\sqrt{t}+\langle Ut\rangle^{\kappa/(\kappa-1)} + {\bf 1}_{M_t\sqrt{t}\ge \langle Ut\rangle \langle x\rangle^{1/\kappa}}
 \left( \frac{M_t\sqrt{t}}{\langle Ut\rangle}\right)^{\kappa} \right)^{
\frac{\alpha}{2}+\frac{1}{\kappa}} \right] \nonumber\\
&\lesssim & K_0 (|x|+\langle Ut\rangle^{\kappa/(\kappa-1)})^{\frac{\alpha}{2}+\frac{1}{\kappa}}. 
\label{eq:bound-u} \EEA
(see proof of Lemma \ref{lem:viscous-admissible} (i)).


\section{Proof of the convergence of the scheme}


The general assumptions on $u_0$ in this main section are:
{\bf \begin{itemize}
\item[(i)] $u_0$ is a $C^2$ function;
\item[(ii)] {\em (a priori bounds on $u_0$, $\nabla u_0$, $\nabla^2 u_0$)} there exist constants $\alpha,\beta\ge 0$ such that, for all $x\in\R^d$,
\BEQ |u_0(x)|\le K_0(1+|x|)^{\frac{\alpha}{2}+\frac{1}{\kappa}}, 
\qquad |\nabla u_0(x)|\le  K_1 (1+|x|)^{\alpha+\frac{2}{\kappa}}, \qquad 
|\nabla^2 u_0(x)|\le K_2 (1+|x|)^{\frac{3}{2}(\frac{\alpha}{2}+\frac{1}{\kappa})} 
\label{eq:a-priori-bound} \EEQ
with $K_0\le U^{\frac{\beta}{2}+1}, K_0\le K_1^{1/2}$, $U\le K_1\le K_2^{2/3}$;
\item[(iii)] $u_0$ coincides outside the union of annuli $\cup_{i\ge 1} {\cal A}_i$ with
an initial velocity $\tilde{u}_0$ satisfying (Hyp1),
\end{itemize} }
\noindent annuli $({\cal A}_i)_{i\ge 1}$ being as in Definition \ref{def:dangerous-zones}.

\medskip\noindent Note that this set of assumptions is precisely that of Theorem \ref{th:viscous-admissible}, plus
some extra a priori bounds on $\nabla u_0$, $\nabla^2 u_0$. We let $M_t:=1+\frac{\sup_{0\le s\le t} |B_s|}{\sqrt{t}}$ as in the previous sections. Generalizing (iii), we may assume
that the sequence of random characteristics $(Y^{(m)}(t;\cdot,x))_{m\ge 0}$ satisfies some weaker form
of the conclusions of Theorem \ref{th:viscous-admissible},

{\bf \begin{itemize}
\item[(iii)']  random characteristics $(Y^{(m)}(\cdot;\cdot,\cdot))_{m\ge 0}$ 
obey the following estimates,
\end{itemize} 

 \BEQ |Y^{(m)}(t;s,x)-x|\le  (C_{\kappa}-1) \langle Ut\rangle  \max(\langle Ut\rangle^{\kappa/(\kappa-1)},|x|)^{1/\kappa} \label{hyp:1}\EEQ
 \qquad\qquad\qquad\qquad\qquad\qquad if 
$M_t\sqrt{t}\le \max(\langle Ut\rangle ^{\kappa/(\kappa-1)}, \langle Ut \rangle \langle x\rangle^{1/\kappa});$

\BEQ |Y^{(m)}(t;s,x)-x|\lesssim (M_t\sqrt{t})^{\kappa'} \label{hyp:2} \EEQ
\qquad\qquad\qquad\qquad\qquad\qquad if
$M_t\sqrt{t}\ge \max(\langle Ut \rangle^{\kappa/(\kappa-1)}, \langle Ut \rangle \langle x\rangle^{1/\kappa})$,

\noindent for some large enough constant $C_{\kappa}>1$, and some exponent $\kappa'\ge 1$ possibly differing from $\kappa$,}

\medskip

\noindent hypothesis (iii) or more generally (iii)'  implying in turn a uniform in $m$  bound on $u^{(m)}$,
\BEQ |u^{(m)}(t,x)|\lesssim K_0(|x|+\langle Ut\rangle^{\kappa/(\kappa-1)})^{\frac{\alpha}{2}+\frac{1}{\kappa}}
\label{hyp:3} \EEQ
(see (\ref{eq:bound-u})), which completes the proof of Theorem 1 in the
Introduction.

\medskip
\noindent We now proceed to prove by induction the bounds on $\nabla u^{(m)},\nabla^2 u^{(m)}, v^{(m)},\nabla v^{(m)}$ collected in Theorem 2 (see section 1.2). All subsequent computations rely exclusively on Feynman-Kac's formula,
Schauder estimates, hypotheses (i),(ii), the bounds on the characteristics, (\ref{hyp:1},\ref{hyp:2}), and their immediate corollary (\ref{hyp:3}).


\subsection{Scheme of proof}


We first want to bound the gradient functions $\nabla u^{(m)}$, $m\ge 0$. 
By using the Feynman-Kac representation and the bounds on the characteristics (\ref{hyp:1},
\ref{hyp:2}), it is easy in the
non-viscous case to
derive local a priori bounds for the gradient in some initial regime $t\le T_{min}(x)$; however,
since $T_{min}(x)\to 0$ when $|x|\to\infty$, one cannot draw from this fact alone any conclusion about  global-in-space, local-in-time regularity of the solution. This works also fine in the viscous case {\em provided} $\alpha=0$, i.e. $u_0$ is {\em sublinear} (or, in other words, {\em if}  (Hyp1) {\em is verified}), and
$\nabla u_0$ {\em subquadratic}, because large deviation estimates (i.e. Gaussian
bounds) for Brownian motion suffice to control the gradient for $t\le T_{min}(x)$. In the
latter case, parabolic Schauder estimates (requiring a non-zero viscosity) make it possible to extend these bounds to
arbitrarily large time. To deal with the general (viscous) case,  we replace eq. (\ref{eq:um}) for
$u^{(m)}$ by a family $u^{(m,n)}$ of {\em penalized}  transport equations, meant as a
smoothened substitute of  the original equation solved on 
dyadic balls $B(0,2^n), n\ge 0$  with Dirichlet boundary conditions. Gradient bounds for the
solutions $u^{(m,n)}$ are easily obtained in some $n$-dependent initial regime $t\le T_n(x)$, and again extended to later times
thanks to Schauder estimates. Then we prove that the series $\sum_n |u^{(m,n)}-u^{(m,n-1)}|$ converges. The same techniques can be repeated to bound second derivatives $\nabla^2 u^{(m)}$  (see \S 4.2). 

In turn we use the uniform estimates for $\nabla u^{(m)}$ found in \S 4.2, together
with those for $u^{(m)}$ (see (\ref{hyp:3})) to bound $v^{(m)}:=
u^{(m)}-u^{(m-1)}$ by simple time integration. For fixed $x$,  we obtain $v^{(m)}(t,x)=O\left(\left(K_1\frac{t}{m} \right)^m\right)$ for $t=O(m/K_1)$ (called: {\em short-time regime}), $O(1)$ otherwise. Thus for fixed $t,x$,
the series $\sum_m |v^{(m)}|$ converges locally uniformly (see \S 4.3).

Finally, repeating the techniques of \S 4.2, we bound $\nabla v^{(m)}$ 
and deduce that the series $\sum_m |\nabla v^{(m)}|$ converges locally uniformly (see \S 4.4). Thus the limit
of the series is a solution of Burgers' equation. 

Note that, by  a standard argument using Schauder's estimates, the solution may be proved to be smooth
for $t>0$. If higher order derivatives of $u_0$ are polynomially bounded, then the regularity
may be proved along the same lines to extend downto $t=0$. In particular, the solution is classical if $u_0$ is $C^2$.


\subsection{Gradient bounds}


We prove in this section the bounds (\ref{intro:nabla-u}), (\ref{intro:nabla2-u})
on $\nabla u^{(m)}$ and $\nabla^2 u^{(m)}$.


\subsubsection{Gradient bounds in the initial regime}


By taking the gradient of (\ref{eq:um}), we get
\BEQ (\partial_t-\Del+u^{(m-1)}\cdot\nabla+\nabla u^{(m-1)})\nabla u^{(m)}=0.
\label{eq:eq-nabla-um} \EEQ
Note that $(\nabla u^{(m-1)})$ is a matrix with entries $(\nabla u^{(m-1)})_{ij}:=
\partial_i u_j^{(m-1)}$. Feynman-Kac formula implies the following representation of
the solution,
\BEQ \nabla u^{(m)}(t,x)=\esper\left[ T\left( e^{-\int_0^t \nabla u^{(m-1)}(t-s,X^{(m)}(t;s,x))\, ds} \right) \nabla u_0(X^{(m)}(t,x)) \right] \label{eq:nabla-um} \EEQ
where $T(\cdot)$ is the time-ordering operator, namely,
\BEQ T\left( e^{\int_0^t B(s) ds} \right):=\sum_{n\ge 0} \int_{t>s_1>\ldots>s_n>0}
B(s_1)\ldots B(s_n) ds_1\ldots ds_n \EEQ
is the solution at time $t$ of the matrix-valued ode $\frac{d}{dt}M(t)=B(t)M(t)$ started
from the identity. We will be happy with the simple bound in terms of matrix norm
$||\cdot||$, $||M(t)||\le \exp\left(\int_0^t ||B(s)||ds\right).$

Let us illustrate this for $m=0,1$. First
\BEQ \nabla u^{(0)}(t,x)=\esper[\nabla u_0(x+B_t)]=e^{t\Del}\nabla u_0(x) \EEQ
hence (see (\ref{eq:etDely1/kappa}),(\ref{eq:U1/2})) 
\BEQ |\nabla u^{(0)}(t,x)|\lesssim K_1 (1+\sqrt{t}+|x|)^{\alpha+\frac{2}{\kappa}} \lesssim
K_1 (|x|+\langle Ut\rangle^{\kappa/(\kappa-1)})^{\alpha+\frac{2}{\kappa}}.\EEQ
As in \S 2.3, this bound may also be found directly without using the explicit solution
for $u^{(0)}$; namely,
\BEQ |\nabla u^{(0)}(t,x)|\le K_1 \esper[(1+|x|+|B_t|)^{\alpha+\frac{2}{\kappa}}]\lesssim K_1 (1+\sqrt{t}+|x|)^{\alpha+\frac{2}{\kappa}}.\EEQ

Next, we consider the case $m=1$. At this stage one readily understands that the
representation (\ref{eq:nabla-um}) alone does not allow an inductive bound, uniform in $m$,
of $\nabla u^{(m)}(t,x)$ for $t\le T(x)$, where $T(x)>0$ is any {\em deterministic} (possibly $x$-dependent)
time. Namely, assuming $\kappa'\ge \kappa$ to make a case, the function in the time-ordered exponential scales for $m=1$ roughly like 
\BEQ tK_1 \left(|x|+ \langle Ut\rangle^{\kappa/(\kappa-1)} + (M_t\sqrt{t})^{\kappa} \right)^{\alpha+\frac{2}{\kappa}} \gtrsim F(t,M_t):=t
K_1 (M_t\sqrt{t})^{2+\kappa\alpha} \EEQ
 for  $t$ small, i.e.  $Ut\le 1$, and  $M_t$ large, i.e. $M_t\sqrt{t}\ge \max(\langle Ut\rangle ^{\kappa/(\kappa-1)}, \langle Ut \rangle \langle x\rangle^{1/\kappa})\approx (1+|x|)^{1/\kappa}$. Hence $F(t,M_t)$ grows for fixed $t$ roughly like $M_t^{\gamma}$, with
$\gamma=2+\kappa\alpha>2$ as soon as $\alpha>0$, which gives
seemingly an infinite average for the
exponential factor   (compare with Gaussian queue (\ref{eq:Mt})). On the other hand (see more details below), we note that in the 'normal' regime
where (assuming $Ut\le 1$) $Y^{(1)}(t;s,x)\in B_{\kappa,CU}(t,x)$, $|Y^{(1)}(t;s,x)|\lesssim \langle x\rangle$, the function in the exponential scales roughly like $tK_1 (1+|x|)^{\alpha+\frac{2}{\kappa}}$. By reference to this case we let, with $C>1$ large enough

\begin{Definition}
 $T_{min}(x):=\left(C^3 K_1 (1+|x|)^{\alpha+\frac{2}{\kappa}}\right)^{-1}$ 
 \label{def:Tmin}
\end{Definition}

\noindent and terminate this somewhat sloppy discussion by some detailed computations.

\medskip \noindent
Let $t\le T_{min}(x)$ (implying in particular $Ut\le 1$ by Hypothesis (ii)), and $\Omega:={\bf 1}_{M_{T_{min}(x)}\ge 
\sqrt{K_1} (1+|x|)^{\frac{\alpha}{2}+\frac{2}{\kappa}} }$. On $\Omega^c$ one has
$M_t\sqrt{t}\le M_{T_{min}(x)}\sqrt{T_{min}(x)}\le (1+|x|)^{1/\kappa}\approx \max(\langle UT_{min}(x)\rangle^{\kappa/(\kappa-1)},
\langle UT_{min}(x)\rangle \langle x\rangle^{1/\kappa})$, hence one is in the 'normal' regime where convection
dominates over diffusion. Then $|Y^{(1)}(t;s,x)-x| \lesssim \langle Ut\rangle\  \max(\langle Ut\rangle^{\kappa/(\kappa-1)},|x|)^{1/\kappa}\lesssim (1+|x|)^{1/\kappa}$, 
hence $|Y^{(1)}(t;s,x)|,|X^{(1)}(t;s,x)|\lesssim 1+|x|$ as pointed out earlier,
 and 
 \BEQ \int_0^t |\nabla u^{(0)}(t-s,X^{(1)}(t;s,x))|\, ds \lesssim  tK_1 (1+|x|)^{\alpha+\frac{2}{\kappa}}\lesssim  T_{min}(x) K_1(1+|x|)^{\alpha+\frac{2}{\kappa}} \le 1 \label{eq:Tmin} \EEQ
  for $C$ large enough, as required. Similarly, $|\nabla u_0(X^{(1)}(t,x))|\lesssim
 K_1 (1+|x|)^{\alpha+\frac{2}{\kappa}}$. On the whole we have proved:
\BEA I^c &:=& \Big| \esper\left[ {\bf 1}_{\Omega^c} T\left( e^{-\int_0^t \nabla u^{(0)}(t-s,X^{(1)}(t;s,x))\, ds} \right) \nabla u_0(X^{(1)}(t,x)) \right] \Big| \nonumber\\
&\lesssim &   K_1 (1+|x|)^{\alpha+\frac{2}{\kappa}},\EEA
a bound comparable to the a priori bound (\ref{eq:a-priori-bound}) for $u_0$.

However for the time being, we fall short of proving a bound for $|\nabla u^{(1)}(t,x)|$
for $t\le T_{min}(x)$ since we have disregarded the event $\Omega$. The reason is that
we have not used the {\em regularizing effect of diffusion}.

\bigskip

We henceforth develop a more comprehensive strategy of proof, incorporating parabolic
Schauder estimates.

\medskip

{\em By induction we assume that for some large enough constant $C>1$,}

\medskip

{\bf (Induction hypothesis) }
\BEQ |\nabla u^{(m-1)}(t,x)|\le C^2 K_1 (|x|+\langle Ut\rangle^{\kappa/(\kappa-1)})^{\alpha+\frac{2}{\kappa}}. \label{eq:induction-nabla-u} \EEQ

\bigskip
{\em The constant $C$ in (\ref{eq:induction-nabla-u}) is the same as in the definition of
$T_{min}(x)$ (see Definition \ref{def:Tmin}), and also the same as that appearing in the
bounds for $\nabla^2 u^{(m)}$ (see (\ref{eq:induction-nabla2-u})), $v^{(m)}$
(see (\ref{eq:ind-v})) and $\nabla v^{(m)}$ (see (\ref{eq:ind-nabla-v})).} It should be
large enough to satisfy various requirements turning up in the course of the proofs. The
important point to be checked carefully is that it may be chosen {\em uniform in $m$}.

\medskip

\noindent We fix some smooth function $\chi:\R_+\to\R_+$ such that $\chi\big|_{[0,1]}=0$ and $\chi\big|_{[2,+\infty)}=1$, and let $\chi^{(n)}(|x|):=\chi(2^{-n}|x|)$, $n\ge 0$.

\begin{Definition} \label{def:wmn}
\begin{itemize}
\item[(i)] For $n\in\N$,
let $u^{(m,n)}:\R_+\times \R^d\to\R^d$ be the solution  of the transport
equation $(\partial_t-\Del+u^{(m-1)}(t,x)\cdot\nabla)u^{(m,n)}(t,x)=  -2C^2 K_1 (2(1+|x|^2))^{\frac{\alpha}{2}+\frac{1}{\kappa}} \chi^{(n)}(|x|)\,  u^{(m,n)}(t,x)$ with initial condition
$u^{(m,n)}(t=0)=u_0$.
\item[(ii)] Let  $\underline{u}^{(m,n)}:=u^{(m,n)}-u^{(m,n-1)}$ ($n\ge 1$).
\end{itemize}
\end{Definition}

Let us write for short $F_n(x):=2C^2 K_1 (2(1+|x|^2))^{\frac{\alpha}{2}+\frac{1}{\kappa}} \chi^{(n)}(|x|)$. The main properties of $F_n$ are the following:
$F_n(x)\ge 0$, $F_n$ is smooth and:
\begin{itemize}
\item[(i)] $F_n(x)=0$ if $|x|\le 2^n$;
\item[(ii)] $F_n(x)\ge 2C^2 K_1(1+|x|)^{\alpha+\frac{2}{\kappa}}\ge 2|\nabla u^{(m-1)}(t,x)|$ if $Ut\le 1$ and  $|x|\ge 2^{n+1}$;
\item[(iii)] $|\nabla F_n(x)|\lesssim C^2 K_1 (1+|x|)^{\alpha+\frac{2}{\kappa}-1} {\bf 1}_{|x|\ge 2^{n}}$, \ \ $|\nabla^2 F_n(x)|\lesssim C^2 K_1 (1+|x|)^{\alpha+\frac{2}{\kappa}-2} {\bf 1}_{|x|\ge 2^{n}}$.
\end{itemize}
As it happens (see below), the dampening of the solution for $|x|$ large is strong enough to ensure a rapid fall-off outside the ball $B(0,2^n)$; compared to more conventional Dirichlet boundary conditions, this has the advantage of avoiding uncontrollable boundary
effects.

The Feynman-Kac representation for $u^{(m,n)}$ is
\BEQ u^{(m,n)}(t,x)=\esper\left[   u_0(X^{(m)}(t,x)) e^{- \int_0^t ds\,
F_n(X^{(m)}(t;s,x))}\right].\EEQ

By subtracting, one gets
\BEQ \underline{u}^{(m,n)}(t,x)=\esper\left[ {\bf 1}_{X^{(m)}(t;\cdot,x)\not\subset B(0,2^{n-1})}
 u_0(X^{(m)}(t,x)) \left( e^{-\int_0^t ds\,
F_n(X^{(m)}(t;s,x))} -  e^{- \int_0^t ds\,
F_{n-1}(X^{(m)}(t;s,x))} \right) \right] \label{eq:w} \EEQ
where $X^{(m)}(t;\cdot,x):=\{ X^{(m)}(t;s,x), \ 0\le s\le t\}$ is the image of the characteristic.

Differentiating, we get
\BEQ (\partial_t-\Del+u^{(m-1)}\cdot\nabla)\nabla u^{(m,n)}(t,x)=-(\nabla u^{(m-1)}(t,x)+F_n(x))
\nabla u^{(m,n)}(t,x)- \nabla F_n(x) u^{(m,n)}(t,x) \EEQ
yielding the  Feynman-Kac representation
\BEQ \nabla \underline{u}^{(m,n)}(t,x) \equiv w_1^{(m,n)}(t,x)-\int_0^t ds\, \left(w_2^{(m,n)}(t;s,x)+w_3^{(m,n)}(t;s,x)+
w_4^{(m,n)}(t;s,x) \right),\EEQ
with (letting $X^{(m)}(t;\le s,x):=\{ X^{(m)}(t;s',x), \ 0\le s'\le s\}$) :
\BEA &&  w_1^{(m,n)}(t,x):=\esper\left[ {\bf 1}_{X^{(m)}(t;\cdot,x)\not\subset B(0,2^{n-1})}
 \left( e^{- \int_0^t ds\,
F_n(X^{(m)}(t;s,x))} -  e^{- \int_0^t ds\,
 F_{n-1}(X^{(m)}(t;s,x))} \right) \right. \nonumber\\
&& \qquad \qquad \qquad \left.  T\left(e^{-\int_0^t ds\,  \nabla u^{(m-1)}(t-s,X^{(m)}(t;s,x))} \right) 
\nabla u_0(X^{(m)}(t,x)) \right]  \label{eq:w1} \EEA
\BEA && w_2^{(m,n)}(t;s,x):=\esper\left[ {\bf 1}_{X^{(m)}(t;\le s,x)\not\subset B(0,2^{n-1})}
 \left( e^{- \int_0^s ds'\,
F_n(X^{(m)}(t;s',x))} -  e^{- \int_0^s ds'\,
 F_{n-1}(X^{(m)}(t;s',x))} \right) \right. \nonumber\\
&& \qquad \qquad \qquad \left.  T\left(e^{-\int_0^s ds'\, \nabla u^{(m-1)}(t-s',X^{(m)}(t;s',x))} \right) \nabla F_{n-1}(X^{(m)}(t;s,x))  u^{(m,n-1)}(t-s,X^{(m)}(t;s,x)) \right];  \nonumber\\  \label{eq:w2}
\EEA 
\BEA &&  w_3^{(m,n)}(t;s,x):= \esper\left[ {\bf 1}_{X^{(m)}(t;\le s,x)\not\subset B(0,2^{n-1})}
  e^{-\int_0^s ds'\,
F_n(X^{(m)}(t;s',x))} \right. \nonumber\\
&& \qquad\qquad \left.   T\left(e^{-\int_0^s ds'\, \nabla u^{(m-1)}(t-s',X^{(m)}(t;s',x))} \right) \nabla(F_n-F_{n-1})(X^{(m)}(t;s,x))   u^{(m,n-1)}(t-s,X^{(m)}(t;s,x))
\right]; \nonumber\\ \label{eq:w3} \EEA 
 \BEA && w_4^{(m,n)}(t;s,x):=\esper\left[  e^{-\int_0^s ds'\,
F_n(X^{(m)}(t;s',x))}  T\left(e^{-\int_0^s ds'\, \nabla u^{(m-1)}(t-s',X^{(m)}(t;s',x))} \right) \right. \nonumber\\
&& \qquad\qquad\qquad\qquad \left.
 \nabla F_n(X^{(m)}(t;s,x)) \underline{u}^{(m,n)}(t-s,X^{(m)}(t;s,x)) \right], \label{eq:w4} \EEA

as deduced from the Feynman-Kac representation for $\nabla u^{(m,n)}(t,x)$,

\BEQ \nabla u^{(m,n)}(t,x) \equiv u_1^{(m,n)}(t,x)-\int_0^t ds\, u_2^{(m,n)}(t;s,x),\EEQ
where 
\BEQ u_1^{(m,n)}(t,x):=\esper\left[ 
  e^{- \int_0^t ds\,
F_n(X^{(m)}(t;s,x))} T\left(e^{-\int_0^t ds\,  \nabla u^{(m-1)}(t-s,X^{(m)}(t;s,x))} \right) 
\nabla u_0(X^{(m)}(t,x)) \right]  \label{eq:u1} \EEQ

\BEA &&  u_2^{(m,n)}(t;s,x):=\esper\left[ 
  e^{- \int_0^s ds'\,
F_n(X^{(m)}(t;s',x))} T\left(e^{-\int_0^s ds'\,  \nabla u^{(m-1)}(t-s',X^{(m)}(t;s',x))} \right)  \right. \nonumber\\
&& \qquad\qquad\qquad \left.
\nabla F_n(X^{(m)}(t;s,x))  u^{(m,n)}(t-s,X^{(m)}(t;s,x)) \right] \label{eq:u2} \EEA

\medskip
We shall now bound: $u^{(m,n)}(t,y)$, $\underline{u}^{(m,n)}(t,y)$ $(y\in\R^d)$ for
$t\le T_{min}(0)=(C^3 K_1)^{-1}$; and  each of the terms contributing to $\nabla \underline{u}^{(m,n)}(t,x)$ for $\langle x\rangle \le (2C_{\kappa})^{-2} 2^{n-1}$ and $t<T_n$,  where
\BEQ T_n:=\left(C^3 K_1 (2^n)^{\alpha+\frac{2}{\kappa}} \right)^{-1}, \qquad n\ge 0.\EEQ
Note that $T_n\approx T_{min}(2^n)$.

The {\em main point} to be understood is that the events $ \left( X^{(m)}(t;\le s,x)\not\subset
B(0,2^{n-1}) \right) $,  figuring inside the expectations defining $\underline{u},w_1,w_2$ and $w_3$, are extremely unlikely for $n$ large. Namely, choose $C_{\kappa}$ large enough; by hypothesis,
\BEQ |X^{(m)}(t;s,x)-x|\le |Y^{(m)}(t;s,x)-x|+M_t\sqrt{t}\le  (C_{\kappa}-1) \langle x\rangle ^{1/\kappa} +O((M_t\sqrt{t})^{\kappa'}) \label{eq:main-point} \EEQ
 for $t\le U^{-1}$ (recall $\kappa'\ge 1$). From this we conclude: if 
$\langle x\rangle \le (2C_{\kappa})^{-1} 2^{n-1}$  (hence in particular, $2^n\ge 4C_{\kappa}\gg 1$), and $t\le T_{min}(0)$,
\BEQ M_t\gtrsim \frac{(2^n)^{\kappa'}}{\sqrt{T_{min}(0)}} \ge C^{3/2} \sqrt{K_1} \  (2^n)^{1/\kappa'}, \label{eq:4.25} \EEQ
an event of probability $O(e^{-cC^3}) O(e^{-cK_1}) O(e^{-c (2^n)^{2/\kappa'}})$.

\begin{itemize}
\item[(i)] (bound for $u^{(m,n)}(t,x)$, $t\le T_{min}(0)$) We replace $u^{(m,n)}(t,x)=\esper[\, \cdot\, ]$ with $\esper[ {\bf 1}_{X^{(m)}(t,x)\in B(0,2^{n-1})} \, \cdot\, ]+ \sum_{p\ge n-1} \esper[ {\bf 1}_{X^{(m)}(t,x)\in B(0,2^{p+1})\setminus
B(0,2^p)} \, \cdot\, ]$. Since $|u_0(X^{(m)}(t,x))|\le K_0 (1+|X^{(m)}(t,x)|)^{\frac{\alpha}{2}+\frac{1}{\kappa}}$, the first and main term is a $O(K_0 (2^n)^{\frac{\alpha}{2}+\frac{1}{\kappa}})$. Subsequent terms are 
$\lesssim K_0 (2^p)^{\frac{\alpha}{2}+\frac{1}{\kappa}} \ \cdot\  O(e^{-cK_1}) O(e^{-c (2^p)^{2/\kappa'}}) \lesssim e^{-c' (2^p)^{2/\kappa'}}$,
summing up to $O(1)$. 

Let us also bound $u^{(m,n-1)}(t,y)$, with $y\in\R^d$ (see (\ref{eq:w3})). If $|y|\ll 2^n$, the
bound is $O(K_0 (2^n)^{\frac{\alpha}{2}+\frac{1}{\kappa}})$ as before. Otherwise, by
a similar reasoning as in (\ref{eq:4.25}), {\em the events} $\left( |X^{(m)}(t-s,y)-y|\gg |y|
 \right) $ {\em  are extremely unlikely for $n$ large}, hence we  get a {\em polynomial bound},
$|u^{(m,n-1)}(t,y)|\lesssim K_0 (1+|y|)^{\frac{\alpha}{2}+\frac{1}{\kappa}}$.

\item[(ii)] (bound for $\underline{u}^{(m,n)}(t,x)$, $t\le T_{min}(0)$) The exponentially small factors in the right-hand side of
(\ref{eq:w}) are not needed for the bound.  We replace $\underline{u}^{(m,n)}(t,x)=\esper[\, \cdot\, ]$ with $\esper[ {\bf 1}_{X^{(m)}(t,x)\in B(0,2^n)} \, \cdot\, ]$ $+ \sum_{p\ge n} \esper[ {\bf 1}_{X^{(m)}(t,x)\in B(0,2^{p+1})\setminus
B(0,2^p)} \, \cdot\, ]$.  Since $|u_0(X^{(m)}(t,x))|\le K_0 (1+|X^{(m)}(t,x)|)^{\frac{\alpha}{2}+\frac{1}{\kappa}}$, the first and main term is a $O(e^{-c' (2^n)^{2/\kappa'}})$. Subsequent terms are $O(e^{-c' (2^p)^{2/\kappa'}})$, summing up also
to $O(e^{-c' (2^n)^{2/\kappa'}})$.

Let us also bound $\underline{u}^{(m,n)}(t-s,y)$ with $y\in\R^d$ (see (\ref{eq:w4})). Reasoning as in (i),
we find: $|\underline{u}^{(m,n)}(t-s,y)|= O(e^{-c' (2^n)^{2/\kappa'}})$ if 
$\langle y\rangle \le (2C_{\kappa})^{-1} 2^{n-1}$, otherwise $|\underline{u}^{(m,n)}(t-s,y)|\lesssim K_0 (1+|y|)^{\frac{\alpha}{2}+\frac{1}{\kappa}}$.

\item[(iii)] (bound for $w_1^{(m,n)}(t,x)$, $t\le T_n$) First we use the matrix bound
$$ ||T\left(e^{-\int_0^t ds\,  \nabla u^{(m-1)}(t-s,X^{(m)}(t;s,x))} \right) 
|| \le \exp \left( \int_0^t ds\,  |\nabla u^{(m-1)}(t-s,X^{(m)}(t;s,x))| \right).$$ Whenever
$|X^{(m)}(t;s,x)|>2^{n+1}$, $e^{-F_{n'}(X^{(m)}(t;s,x)) + |\nabla u^{(m-1)}(t-s,X^{(m)}(t;s,x))|}\le 1$,
$n'=n,n-1$. On the other hand, if $|X^{(m)}(t;s,x)|\le 2^{n+1}$, then $|\nabla u^{(m-1)}(t-s,X^{(m)}(t;s,x))|\le C^2 K_1 (1+2^{n+1})^{\alpha+\frac{2}{\kappa}}$. Thus the product of the exponential factors is $\le \exp\  O\left(T_n\, \cdot\,  C^2 K_1 (1+2^{n+1})^{\alpha+\frac{2}{\kappa}}
\right)\le e$ for $C$ large enough. Then the product of the characteristic function with $\nabla u_0(X^{(m)}(t,x))$ is bounded by $O(e^{-c'(2^n)^{2/\kappa'}})$ by the same arguments as in (ii).

\item[(iv)] (bound for $w_2^{(m,n)}(t,x)$ and $w_3^{(m,n)}(t,x)$, $t\le T_n$) The time-ordered exponential
is compensated as in (iii).  Proceeding as in (ii), we see that the main contribution comes from the case $X^{(m)}(t;s,x)\in B(0,2^n)$.  Then $|\nabla F_{n'}(X^{(m)}(t;s,x)|
\lesssim C^2 K_1 (2^n)^{\alpha+\frac{2}{\kappa}-1}$, while  $|u^{(m,n-1)}(t-s,X^{(m)}(t;s,x))|\lesssim K_0 (2^n)^{\frac{\alpha}{2}+\frac{1}{\kappa}}\le K_1^{1/2}
(2^n)^{\frac{\alpha}{2}+\frac{1}{\kappa}}$. 
Taking the product with  the characteristic function yields $O(e^{-c'(2^n)^{2/\kappa'}})$.

\item[(v)] (bound for $w_4^{(m,n)}(t,x)$, $t\le T_n$)  Replace $w_4^{(m,n)}(t;s,x)=\esper[\cdot]$ with
$\esper\left[{\bf 1}_{X^{(m)}(t;s,x)\not\in B(0,(2C_{\kappa})^{-1} 2^{n-1})}  \ \cdot\  \right]$ $+\esper\left[{\bf 1}_{X^{(m)}(t;s,x)\in  B(0,(2C_{\kappa})^{-1} 2^{n-1})} \ \cdot\ \right]$. The first term
is bounded  by $O(e^{-c' (2^n)^{2/\kappa'}})$ as in (ii), since (by hypothesis) $|x|\le (2C_{\kappa})^{-2} 2^{n-1}$. Assume on the other hand $X^{(m)}(t;s,x)\in B(0,(2C_{\kappa})^{-1}2^{n-1})$;
 then $|\underline{u}^{(m,n)}(t-s,X^{(m)}(t;s,x)|=O(e^{-c'(2^n)^{\alpha+\frac{2}{\kappa}}})$, as proved in (ii).

\end{itemize}

Leaving aside the bounds for $u^{(m,n)}$ and $\underline{u}^{(m,n)}$, which shall be
used in \S 4.2.2 below, we have proved:
\BEQ |\nabla \underline{u}^{(m,n)}(t,x)|\lesssim e^{-c'(2^n)^{2/\kappa'}}, \label{eq:4.29} \EEQ
valid for $t\le T_n$ and $\langle x\rangle \le (2C_{\kappa})^{-2} 2^{n-1}$.

\bigskip

For a {\em given dyadic slice} 
\BEQ x\in B(0,2^p)\setminus B(0,2^{p-1}), \qquad p\ge 1, \EEQ
 one may apply this result for any  $n\ge n':=p+1+ \lceil 2\log_2(2C_{\kappa}) \rceil$.

We now {\em assume} $|x|\ge (2C_{\kappa})^{\kappa/(\kappa-1)}$ (so that $(C_{\kappa}-1)
\langle x\rangle ^{1/\kappa}\le\half |x|$, see (\ref{eq:main-point})) , fix  $n'':=p-1-\lceil \log_2(2C_{\kappa})\rceil\ge 0$  and write
\BEQ \nabla u^{(m)}(t,x)=\nabla u^{(m,n'')}(t,x)+  (\nabla \underline{u}^{(m,n''+1)}(t,x)+\ldots+
\nabla \underline{u}^{(m,n'-1)}(t,x)) + \sum_{n\ge n'} \nabla \underline{u}^{(m,n)}(t,x) \label{eq:w-decomposition}\EEQ
as a sum of three contributions, in which $|x|$ is {\em large} (first term, $\nabla u^{(m,n'')}(t,x)$), {\em small} (last term, $\sum_{n\ge n'}\nabla \underline{u}^{(m,n)}(t,x)$), or
of the same order as $2^n$.

\medskip

{\em Our purpose is to show that: $|\nabla u^{(m,n)}(t,x)|$ ($n=n''$), 
$|\nabla \underline{u}^{(m,n)}(t,x)|$ $(n=n''+1,\ldots,n'-1)$ are $\lesssim K_1
(1+|x|)^{\alpha+\frac{2}{\kappa}}$, while the remaining terms, $\nabla \underline{u}^{(m,n)}(t,x)$, $n\ge n'$ are negligible} (see (\ref{eq:4.29})). The problem, however,
is that, for the time being, we shall be able to prove these only for $t\le T_n$.
Since $T_n\to_{n\to\infty} 0$, we cannot say anything about the sum in (\ref{eq:w-decomposition}) till we extend these bounds to arbitrary time (see next subsection).

\medskip

Consider now the first term $(x$ {\em large}) with $t\le T_{n''}$. As in (i), the contribution coming from the case
$\sup_{0\le s\le t} |X^{(m)}(t;s,x)-x|\ge \frac{2}{3}|x| $ is $O(e^{-c(2^{n''})^{2/\kappa'}})$. In the contrary case,  $|X^{(m)}(t;s,x)|\le 2|x|$ for all $0\le s\le t$, so
  $|\nabla u_0(X^{(m)}(t;s,x))|\lesssim K_1 (1+|x|)^{\alpha+\frac{2}{\kappa}}$, while  
\BEA && \int_0^t ds\, |\nabla F_{n''}(X^{(m)}(t;s,x)) u^{(m,n'')}(t-s,X^{(m)}(t;s,x))|
 \nonumber\\
 && \qquad \qquad\qquad \lesssim T_{n''}\ \cdot\ C^2 K_1 (1+|x|)^{\alpha+\frac{2}{\kappa}-1} \ \cdot\ K_0 (1+|x|)^{\frac{\alpha}{2}+\frac{1}{\kappa}} \lesssim
 K_0 (1+|x|)^{\frac{\alpha}{2}+\frac{1}{\kappa}-1}.\EEA
 All together we have found: $|\nabla u^{(m,n'')}(t,x)|\lesssim K_1 (1+|x|)^{\alpha+\frac{2}{\kappa}}$.

\medskip

Consider finally the  finite number of terms $n=n''+1,\ldots,n'-1$ for which $|x|\approx 2^n$. 
Reasoning as in (i) we may assume that $|X^{(m)}(t,x)|, |X^{(m)}(t;s,x)| \lesssim |x|$
in the above formulas, whence 
\BEQ w_1^{(m,n)}(t,x)\lesssim K_1 (1+|x|)^{\alpha+\frac{2}{\kappa}}; \EEQ
\BEQ w_2^{(m,n)}(t;s,x),w_3^{(m,n)}(t;s,x),w_4^{(m,n)}(t;s,x)\lesssim CK_1 (1+|x|)^{\alpha+\frac{2}{\kappa}-1} \ \cdot\ K_0 (1+|x|)^{\frac{\alpha}{2}+\frac{1}{\kappa}} \EEQ
and finally,
\BEQ |\nabla \underline{u}^{(m,n)}(t,x)|\lesssim K_1 (1+|x|)^{\alpha+\frac{2}{\kappa}} + T_n \  CK_1 (1+|x|)^{\alpha+\frac{2}{\kappa}-1} \ \cdot\ K_0 (1+|x|)^{\frac{\alpha}{2}+\frac{1}{\kappa}}
 \lesssim K_1 (1+|x|)^{\alpha+\frac{2}{\kappa}}. \EEQ
Clearly the estimates are the same as for $\nabla u^{(m,n'')}$, so in the sequel we
shall group together these two terms and rewrite (\ref{eq:w-decomposition}) as
\BEQ \nabla u^{(m)}(t,x)=\nabla u^{(m,n'-1)}(t,x)+  \sum_{n\ge n'} \nabla \underline{u}^{(m,n)}(t,x) \label{eq:w-decomposition-bis}\EEQ
Note that 
 all these arguments are easily adapted to the case $|x|<(2C_{\kappa})^{\kappa/(\kappa-1)}$ provided $C$ is large enough (take $n''=0$).
 
\bigskip 
 
Let us recapitulate. Summing the three contributions from (\ref{eq:w-decomposition}), or the two
contributions from (\ref{eq:w-decomposition-bis}), we see
that  (again, provided $C$ is large enough) our induction hypothesis (\ref{eq:induction-nabla-u}) should hold at rank $m$,
except that our gradient bounds should be proven to hold for all $t>0$; and to start with, if possible, for all $t$ less than some uniform stopping time, $t\le T_{min}(0)=(C^3 K_1)^{-1}$. This is
precisely what we do in the next paragraph.


\subsubsection{Large-time bounds for the gradient}


For $t$ away from the time origin, bounds for the gradient rest on Schauder estimates. We
use a quantitative form of these proved by us in \cite{Unt-Bur1}. Let us quote the
result for the sake of the reader. More detailed bounds are proved in \cite{Unt-Bur1},
 Proposition 4.5.

\begin{Proposition} \cite{Unt-Bur1} \label{prop:Schauder}
Let $v$ solve the linear parabolic PDE 
\BEQ (\partial_t-\Del+a(t,x))u(t,x)=b(t,x)\cdot
\nabla u(t,x) + f(t,x) \EEQ
 on the "parabolic ball" 
$Q^{(j)}=Q^{(j)}(t_0,x_0):=\{(t,x)\in\R\times\R^d;\ t_0-2^{j}\le t\le t_0, x\in \bar{B}(x_0,2^{j/2})\}$. If $u$ is bounded, $a\ge 0$,  
\BEQ ||f||_{\gamma,Q^{(j)}}:=\sup_{(t,x),(t',x')\in Q^{(j)}} \frac{|f(t,x)-f(t',x')|}{|x-x'|^{\gamma}+|t-t'|^{\gamma/2}}<\infty \EEQ for some $\gamma\in(0,1)$,
and similarly $||a||_{\gamma,Q^{(j)}}, ||b||_{\gamma,Q^{(j)}}<\infty$,
 then 
\BEQ 
\sup_{Q^{(j-1)}} |\nabla u| \lesssim 2^{j/2} R_b^{-1} \left( 
 2^{j\gamma/2} ||f||_{\gamma,Q^{(j)}}+(2^{j\gamma} R_b^{-1} ||b||_{\gamma,Q^{(j)}}^2
 + 2^{j\gamma/2} ||a||_{\gamma,Q^{(j)}} + 2^{-j}) \sup_{Q^{(j)}} |u| \right), \label{eq:S1} \EEQ 
 
 \BEQ 
\sup_{Q^{(j-1)}} |\partial_t u|, \sup_{Q^{(j-1)}} |\nabla^2 u| \lesssim  R_b^{-1} \left( 
 2^{j\gamma/2}  ||f||_{\gamma,Q^{(j)}}+(  2^{j\gamma} R_b^{-1} ||b||_{\gamma,Q^{(j)}}^2+ 2^{j\gamma/2} ||a||_{\gamma,Q^{(j)}}+2^{-j}) \sup_{Q^{(j)}} |u| \right), 
\label{eq:S2} \EEQ
 where $R_b:=\left(1+2^{j/2}|b(t_0,x_0)|\right)^{-1}$.
\end{Proposition}

\medskip

Fix $\gamma\in(0,1)$ and  $x\in B(0,2^p)\setminus B(0,2^{p-1})$, $p\ge 1$ in a given
dyadic slice. Define $n':=p+1+\lceil 2\log_2 (2C_{\kappa})\rceil$ as in (\ref{eq:w-decomposition})). Recall we have shown: $|\nabla u^{(m,n'-1)}(t,x)|\lesssim K_1(1+|x|)^{\alpha+\frac{2}{\kappa}}$ for $t\le T_{n'-1}$, and $|\nabla \underline{u}^{(m,n)}(t,x)|\lesssim e^{-c'(2^n)^{2/\kappa'}}$ $(n\ge n')$ for $t\le T_n$.

\begin{enumerate}
\item We consider first the initial regime $t\le T_{min}(0)$, where bounds (i),(ii) 
for $u^{(m,n)}, \underline{u}^{(m,n)}$ hold (see \S 4.2.1). 
Decomposing $u$ as $u^{(m,n'-1)}+ \sum_{n\ge  n'} \underline{u}^{(m,n)}$,  we apply  Proposition \ref{prop:Schauder},  (i) to $u^{(m,n'-1)}$ on $Q:=Q^{(\log_2 T_{n'-1})}(t,x)$, $t\ge T_{n'-1}$ ($x$ {\em large}); (ii) to $\underline{u}^{(m,n)}$ on $Q:=Q^{(\log_2 T_n)}(t,x)$,
$t\ge  T_n$  for $n\ge n'$
($x$ {\em small}),
   with $b:=-u^{(m-1)}$,
$f\equiv 0$
and (i) $a(t,x):=F_n(x)$,
(ii) $a\equiv 0$.

We concentrate on case (i), where $2^j=T_{n'-1}\approx T_p\approx (C^3 K_1 \langle x\rangle^{\alpha+\frac{2}{\kappa}})^{-1}$. Then $R_b^{-1}=1+\sqrt{T_{n'-1}} |u^{(m-1)}(t,x)|\lesssim 1$. By H\"older interpolation, 
\BEA  ||u^{(m-1)}||_{\gamma,Q} &\lesssim &   \left(\sup_Q |u^{(m-1)}|\right)^{1-\gamma} \left( \sup_Q
|\nabla u^{(m-1)}|\right)^{\gamma}\lesssim
 \left( K_0\langle x\rangle^{\frac{\alpha}{2}+\frac{1}{\kappa}} \right)^{1-\gamma}
  \left( C^2 K_1 \langle x\rangle^{\alpha+\frac{2}{\kappa}} \right)^{\gamma}
  \nonumber\\
  &\lesssim& 
  C^{2\gamma}
  K_1^{(1+\gamma)/2} \langle x\rangle^{(\frac{\alpha}{2}+\frac{1}{\kappa})(1+\gamma)}
  \label{eq:4.35} \EEA
since $K_0\le K_1^{1/2}$, and $2^{j\gamma/2} ||a||_{\gamma,Q}\lesssim  2^{j\gamma/2}  \cdot\ C^2 K_1 \langle x\rangle^{\alpha+\frac{2}{\kappa}-\gamma}\lesssim
C^2 K_1 \langle x\rangle^{\alpha+\frac{2}{\kappa}}$, 
  $2^{j\gamma} ||u^{(m-1)}||^2_{\gamma,Q}+2^{j\gamma/2} ||a||_{\gamma,Q^{(j)}}+2^{-j}\lesssim C^3 K_1 \langle x\rangle^{\alpha+\frac{2}{\kappa}}$, $2^{j/2} \sup_Q |u^{(m,n'-1)}|\lesssim C^{-3/2}$, hence Proposition
\ref{prop:Schauder} yields $|\nabla u^{(m,n'-1)}(t,x)|\lesssim C^{3/2} K_1 \langle x\rangle^{\alpha+\frac{2}{\kappa}}$. 

\medskip

Consider now briefly (ii) ($x$ small). Then  one still has
$R_b^{-1}\lesssim 1$, $||u^{(m-1)}||_{\gamma,Q}\lesssim
  C^{2\gamma} K_1^{(1+\gamma)/2} \langle x\rangle^{(\frac{\alpha}{2}+\frac{1}{\kappa})(1+\gamma)}$, while now $T_{n-1}^{1/2} \sup_Q |\underline{u}^{(m,n)}|=O(e^{-c'(2^n)^{2/\kappa'}})$ is
  exponentially small.
  
  \medskip

Summing the two contributions, we see that we have proved what we wanted if $C$ is
large enough: $|\nabla u^{(m)}(t,x)|\le 
C^2 K_1(1+|x|)^{\alpha+\frac{2}{\kappa}}$, for all $t\le T_{min}(0)$ this time.

\item Let now $t\ge T_{min}(0)$. Define 
\BEQ \langle x\rangle_t:=|x|+\langle Ut\rangle^{\kappa/(\kappa-1)}, \qquad 
T_{min}(t,x):=\left(C^3 K_1 \langle x\rangle_t^{\alpha+\frac{2}{\kappa}}\right)^{-1}.
\label{eq:Tmint}  \EEQ
  Apply Proposition \ref{prop:Schauder} directly to $u$
on $Q:=Q^{(\log_2 T_{min}(t,x))}(t,x)$. Then $R_b^{-1}=1+$ \\ $+\sqrt{T_{min}(t,x)}\,  |u^{(m-1)}(t,x)|\lesssim 1$. Instead of (\ref{eq:4.35}) one gets: $||u^{(m-1)}||_{\gamma,Q}\lesssim
C^{2\gamma} K_1^{(1+\gamma)/2} \langle x\rangle_t^{(\frac{\alpha}{2}+\frac{1}{\kappa})(1+\gamma)}$, whence $T_{min}(t,x)^{\gamma} ||u^{(m-1)}||^2_{\gamma,Q} + T_{min}(t,x)^{-1}\lesssim C^3 K_1 
\langle x\rangle_t^{\alpha+\frac{2}{\kappa}}$. Finally, $T_{min}(t,x)^{1/2} \sup_Q |u|\lesssim C^{-3/2}$. Hence Proposition \ref{prop:Schauder} yields for $C$ large enough: $|\nabla u^{(m,n'')}(t,x)|\lesssim C^{3/2} K_1 \langle x\rangle_t^{\alpha+\frac{2}{\kappa}}$.

\end{enumerate}


\subsubsection{Bounds for $\nabla^2 u^{(m)}$}


Unfortunately, in order to prove the convergence of the scheme, we also
need to prove bounds for {\em second-order derivatives} of $u^{(m)}$. However, the
proof proceeds exactly as for the gradient, and we shall only sketch it very roughly.
We want to prove (\ref{intro:nabla2-u}) :

\medskip

{\bf (Induction hypothesis) }
\BEQ |\nabla^2 u^{(m-1)}(t,x)|\le C^4 K_2 (|x|+\langle Ut\rangle^{\kappa/(\kappa-1)})^{3(\frac{\alpha}{2}+\frac{1}{\kappa})}. \label{eq:induction-nabla2-u} \EEQ

\bigskip
Comparing with (\ref{eq:induction-nabla-u}), we see that $|\nabla^2 u^{(m-1)}|$ scales
roughly like $|\nabla u^{(m-1)}|^{3/2}$. This is coherent with the hypothesis
$K_2\ge K_1^{3/2}$. Differentiating once more the equation for $u^{(m,n)}$ (see
Definition \ref{def:wmn}), we get

\BEA && (\partial_t-\Del +(2\nabla u^{(m-1)}(t,x)+F_n(x)) + u^{(m-1)}\cdot\nabla)\nabla^2 u^{(m,n)}(t,x)=
- \nabla(\nabla F_n(x) u^{(m,n)}(t,x)) \nonumber\\
&& \qquad\qquad\qquad -\nabla(\nabla u^{(m-1)}(t,x)+F_n(x))\, \nabla u^{(m,n)}(t,x). 
\label{eq:4.43} \EEA

The Feynman-Kac representation for $\nabla^2 u^{(m,n'-1)}$ or $\nabla^2 \underline{u}^{(m,n)}$, $n\ge n'$ is very much alike that of $\nabla u^{(m,n)}$ or $\nabla \underline{u}^{(m,n)}$, except that there is one more gradient, and there appear supplementary terms due to the last
term in (\ref{eq:4.43}). The exponential multiplicative factor is (up to the coefficient
2 in (\ref{eq:4.43})) the
same as in the case of $\nabla u$, hence may be essentially neglected for
$t<T_n$. Similarly, the convection term may be essentially neglected since
$|X^{(m-1)}(t,x) |\lesssim \langle x\rangle$ with high probability when $t\le T_{min}(0)$.
Thus (considering only the main contribution), for $t\lesssim T_{min}(x)\approx (C^3 K_1 \langle x\rangle^{\alpha+\frac{2}{\kappa}})^{-1}$, and $n=n'-1=\log_2 \langle x\rangle+O(1)$, 
\BEA && |\nabla^2 u^{(m,n)}(t,x)|\lesssim  \sup_{\langle x'\rangle\approx \langle x\rangle} |\nabla^2 u_0(x')| + T_{min}(x)  \ \cdot\  \sup_{0\le t'\le t, \langle x'\rangle  \approx
\langle x\rangle} \left\{  |\nabla(\nabla F_n(x') u^{(m,n)}(t',x'))| \right. \nonumber\\
&& \qquad\qquad\qquad \left. 
+ |\nabla(\nabla u^{(m-1)}(t',x')+F_n(x'))| \  \cdot\ |\nabla u^{(m,n)}(t',x')| 
\right\} .  \label{eq:4.44} \EEA

 In this expression  $|u^{(m,n)}(t',x')|\lesssim CK_0 \langle x\rangle^{\frac{\alpha}{2}+\frac{1}{\kappa}}$, $|\nabla u^{(m,n)}(t',x')|\lesssim
 C^2 K_1 \langle x\rangle^{\alpha+\frac{2}{\kappa}}$, and (by induction)
 $|\nabla^2 u^{(m-1)}(t',x')|\lesssim C^4 K_2 \langle x\rangle^{3(\frac{\alpha}{2}+\frac{1}{\kappa})}$.  The largest terms are obtained by letting the gradient act on $u^{(m,n)}$ since
$|\nabla^2 F_n(x')|\lesssim |\nabla F_n(x')|\lesssim F_n(x')\lesssim C^2 K_1 \langle x\rangle^{\alpha+\frac{2}{\kappa}}$, while bounds on $u^{(m,n)}$ get worse and worse
each time one applies a gradient.  Hence: 
\BEA  |\nabla^2 u^{(m,n)}(t,x)| &\lesssim& K_2 \langle x\rangle^{3(\frac{\alpha}{2}+\frac{1}{\kappa})} \ +\  (C^3 K_1 \langle x\rangle^{\alpha+\frac{2}{\kappa}})^{-1}  \left\{ 
\left( C^2 K_1 \langle x\rangle^{\alpha+\frac{2}{\kappa}} \right)^2
 \ +\ C^4 K_2 \langle x\rangle^{3(\frac{\alpha}{2}+\frac{1}{\kappa})} \ \cdot\ 
 C^2 K_1 \langle x\rangle^{\alpha+\frac{2}{\kappa}}   \right\} \nonumber\\
  &\lesssim & C^3 K_2 \langle x\rangle^{3(\frac{\alpha}{2}+\frac{1}{\kappa})}. \EEA

Taking $C$ large enough one obtains inductively a uniform in $m$ short-time estimate for $\nabla^2 u^{(m,n)}$. For $t$ larger  one must use Schauder estimates as in \S 4.2.2
(see  (\ref{eq:Tmint})). Comparing
(\ref{eq:S1}) with (\ref{eq:S2}) one sees that the bound for $\sup_Q |\nabla^2 u^{(m,n)}|$
or $\sup_Q |\nabla^2 u^{(m)}|$ 
differs from the bound for $\sup_Q |\nabla u^{(m,n)}|$, resp.  $\sup_Q |\nabla u^{(m)}|$ 
only by a multiplicative factor $2^{-j/2}\approx T_p^{-1/2}\approx T_{min}(t,x)^{-1/2}
\approx C^{3/2}  K_1^{1/2} \langle x\rangle_t^{\frac{\alpha}{2}+\frac{1}{\kappa}} \le 
C^{3/2}  \frac{K_2}{K_1} \langle x\rangle_t^{\frac{\alpha}{2}+\frac{1}{\kappa}} $.
 Hence $\sup_Q |\nabla^2 u^{(m)}|\lesssim
C^3 K_2  \langle x\rangle_t^{3(\frac{\alpha}{2}+\frac{1}{\kappa})}$, allowing a bound
uniform in $m$ by induction.


\subsection{Bounds for $v^{(m)}$}


We prove in this section (\ref{intro:v}).
Subtracting eq. (\ref{eq:um}) for $m,m-1$, we find an equation for $v^{(m)}:=u^{(m)}-u^{(m-1)}$,
\BEQ (\partial_t-\Del+u^{(m-1)}(t,x)\cdot\nabla)v^{(m)}(t,x)=f^{(m-1)}(t,x):=-v^{(m-1)}(t,x)\cdot\nabla u^{(m-1)}(t,x). \label{eq:v} \EEQ

\noindent Recall $T_{min}(t,x)=\left(C^3 K_1 \langle x\rangle_t^{\alpha+\frac{2}{\kappa}}\right)^{-1}$ (see (\ref{eq:Tmint})). We assume
\bigskip

\noindent {\bf (Induction hypothesis)} \BEQ |v^{(m-1)}(t,x)|\le CK_0  
(t/(m-1)T_{min}(t,x))^{m-1} (|x|+\langle Ut\rangle^{\kappa/(\kappa-1)})^{\frac{\alpha}{2}+\frac{1}{\kappa}}, \  t>0. \label{eq:ind-v} \EEQ

\medskip

Note that (\ref{eq:ind-v}) is an improvement on (\ref{hyp:3}) only when
$t\le (m-1)T_{min}(t,x)$, i.e. in some initial regime $t\in [0,T_{min}^{(m)}(x)]$, 
where $T_{min}^{(m)}(x)$ is given by an implicit equation (it is easy to show that
$T_{min}^{(m)}(x)\approx (m-1)T_{min}(x)\approx (m-1)
 (C^3 K_1\langle x\rangle^{\alpha+\frac{2}{\kappa}})^{-1}$ for $\langle x\rangle \ge (Ut)^{\kappa/(\kappa-1)}$, in particular for $t\le U^{-1}$,  otherwise
$T_{min}^{(m)}(x)\approx U^{-\lambda} \left(\frac{m-1}{C^3 K_1}\right)^{\mu}$, with
$\lambda=\frac{\kappa\alpha+2}{\kappa(1+\alpha)+1}$, $\mu=\frac{\kappa-1}{\kappa(1+\alpha)+1}<1$). 

\medskip

\noindent Eq. (\ref{eq:v}) also admits a Feynman-Kac representation,
\BEQ v^{(m)}(t,x)=-\int_0^t ds \, \esper\left[ v^{(m-1)}(t-s,X^{(m-1)}(t;s,x)) \cdot 
\nabla u^{(m-1)}(t-s,X^{(m-1)}(t;s,x)) \right].\EEQ

Using the gradient bound, $|\nabla u^{(m-1)}(t,x)|\lesssim C^2 K_1(|x|+\langle Ut\rangle^{\kappa/(\kappa-1)})^{\alpha+\frac{2}{\kappa}}$ and the characteristic
estimate $|X^{(m-1)}(t;s,x)|\lesssim |x|+M_t\sqrt{t}+
\langle Ut\rangle^{\kappa/(\kappa-1)} + {\bf 1}_{M_t\sqrt{t}\ge \langle Ut\rangle 
\langle x\rangle^{1/\kappa}} (M_t\sqrt{t})^{\kappa}$,  we deduce (compare with the
proof of (\ref{eq:bound-u})):
\BEA  |v^{(m)}(t,x)| &\lesssim& \int_0^t ds\, CK_0 ((t-s)/(m-1)T_{min}(t,x))^{m-1}
 \ \cdot \ 
C^2 K_1   (|x|+\langle Ut\rangle^{\kappa/(\kappa-1)})^{3(\frac{\alpha}{2}+\frac{1}{\kappa})}
\nonumber\\
&\lesssim&  K_0  
(t/mT_{min}(t,x))^{m} (|x|+\langle Ut\rangle^{\kappa/(\kappa-1)})^{\alpha+\frac{2}{\kappa}} \nonumber\\
&\le & CK_0  
(t/mT_{min}(t,x))^{m} (|x|+\langle Ut\rangle^{\kappa/(\kappa-1)})^{\alpha+\frac{2}{\kappa}} \EEA
for $C$ large enough.


\subsection{Gradient bounds for $v^{(m)}$}


We prove in this section the bound (\ref{intro:v}) for $\nabla v^{(m)}$.

\medskip

\noindent Differentiating (\ref{eq:v}), one finds
\BEQ (\partial_t-\Del+u^{(m-1)}(t,x)\cdot\nabla + \nabla u^{(m-1)}(t,x))\nabla v^{(m)}(t,x)=\nabla f^{(m-1)}(t,x), \label{eq:nabla-vm} \EEQ
compare with (\ref{eq:eq-nabla-um}), with a right-hand side
\BEQ \nabla f^{(m-1)}(t,x)=\nabla \left( -v^{(m-1)}(t,x)\cdot\nabla u^{(m-1)}(t,x) \right)= -v^{(m-1)}(t,x) \cdot \nabla^2 u^{(m-1)}(t,x)- \nabla v^{(m-1)}(t,x) \cdot  \nabla u^{(m-1)}(t,x).\EEQ

We now proceed as in \S 4.2 to which we refer for the scheme of proof and notations, and define, similarly to Definition \ref{def:wmn}, 

\begin{Definition} 
\begin{itemize}
\item[(i)] For $n\in\N$,
let $v^{(m,n)}:\R_+\times \R^d\to\R^d$ be the solution  of the transport
equation $(\partial_t-\Del+u^{(m-1)}(t,x)\cdot\nabla)v^{(m,n)}(t,x)=  -F_n(x)\,  v^{(m,n)}(t,x)+ f^{(m-1)}(t,x)$ with initial condition
$v^{(m,n)}(0)=0$.
\item[(ii)] Let  $\underline{v}^{(m,n)}:=v^{(m,n)}-v^{(m,n-1)}$ ($n\ge 1$).
\end{itemize}
\end{Definition}

\medskip
Differentiating the equation for $v^{(m,n)}$, we get
\BEQ (\partial_t-\Del+u^{(m-1)}(t,x)\cdot\nabla + (\nabla u^{(m-1)}(t,x)+F_n(x)))\nabla v^{(m,n)}(t,x)=-\nabla F_n(x) v^{(m,n)}(t,x)+ \nabla f^{(m-1)}(t,x). \label{eq:nabla-vmn} \EEQ

The Feynman-Kac representation of $v^{(m,n)}$, $\nabla v^{(m,n)}$, $\underline{v}^{(m,n)}$, $\nabla\underline{v}^{(m,n)}$ are totally similar
to those of $u^{(m,n)}$, $\nabla u^{(m,n)}$, $\underline{u}^{(m,n)}$, $\nabla \underline{u}^{(m,n)}$, with $u^{(m,n-1)}$, $\underline{u}^{(m,n)}$ replaced by their
counterparts $v^{(m,n-1)}$, $\underline{v}^{(m,n)}$ in the expressions for $w_j^{(m,n)}$,
$j=2,3,4$, and the initial condition term $w_1^{(m,n)}(t,x)$ replaced by a contribution due to the right-hand side, $\int_0^t ds\,  w_1^{(m,n)}(t;s,x)$, where
\BEA  && w_1^{(m,n)}(t;s,x):=\esper\left[ {\bf 1}_{X^{(m)}(t;\le s,x)\not\subset B(0,2^{n-1})}
  \left( e^{- \int_0^s ds'\,
F_n(X^{(m)}(t;s',x))} -  e^{- \int_0^s ds'\,
 F_{n-1}(X^{(m)}(t;s',x))} \right) \right. \nonumber\\
&& \qquad \qquad \qquad \left.  T\left(e^{-\int_0^s ds'\, \nabla u^{(m-1)}(t-s',X^{(m)}(t;s',x))} \right) \nabla f^{(m-1)}(t-s,X^{(m)}(t;s,x)) \right]. \EEA

Fix some exponent $\gamma\in(0,1)$, and let $\tilde{T}_{min}(t,x):=\left(C^3 K_2^{2/3}  (|x|+ \langle Ut\rangle^{\kappa/(\kappa-1)})^{\alpha+\frac{2}{\kappa}} \right)^{-1}$. 
 We assume inductively:  

\bigskip

\noindent {\bf (Induction hypothesis)} \ \ \ \ \ \ 
\BEQ |\nabla v^{(m-1)}(t,x)|\le C^3 K_2^{2/3}  
(t/(m-1)\tilde{T}_{min}(t,x))^{\gamma(m-1)/2} (|x|+\langle Ut\rangle^{\kappa/(\kappa-1)})^{\alpha+\frac{2}{\kappa}}, \qquad t\le (m-1)\tilde{T}_{min}(t,x) \label{eq:ind-nabla-v} \EEQ

\medskip

\noindent Let us make two comments at this point. 
First, because $\nabla f^{(m-1)}(t,x)$ involves the second derivative $\nabla^2 u^{(m-1)}$,
which is rougly of order $K_2$ (for $t,x$ small), 
and $K_2^{2/3}\ge K_1$, our bounds are in terms of the larger constant $K_2^{2/3}$ and
not in terms of $K_1$, which also accounts for the replacement of $T_{min}$ by
$\tilde{T}_{min}\le T_{min}$. 
Second, our  bound for $|\nabla v^{(m-1)}(t,\cdot)|$ is in $(t/(m-1))^{\gamma(m-1)/2}$,
$\gamma<1$ 
instead of the naively expected and smaller $(t/(m-1))^{m-1}$ (as found before for $|v^{(m-1)}(t,\cdot)|$) 
for reasons that appear only when applying Schauder estimates (see below).

\medskip

\noindent {\em For $t$ small enough}, bounds for $\nabla v^{(m,n)}$, $\nabla \underline{v}^{(m,n)}$ may
be proved using the Feynman-Kac representation. Let $x\in B(0,2^p)\setminus B(0,2^{p-1})$
$(p\ge 1)$ and $n':=p+1+\lceil 2\log_2 (2C_{\kappa})\rceil$ as in (\ref{eq:w-decomposition}). 
 We refer to the computations in \S 4.2.3. Considering only the main contribution, (\ref{eq:4.44}) is replaced with
\BEA && |\nabla v^{(m,n)}(t,x)|  \lesssim  \int_0^{t} dt'\,  \sup_{\langle
x'\rangle \approx \langle x\rangle} \left(  |\nabla F_n(x')| \  |v^{(m,n)}(t',x')| +
|\nabla f^{(m-1)}(t',x')| \right) \nonumber\\
&&\lesssim  
 \int_0^{t} dt'
 \left\{ C^2 K_2^{2/3}  \langle x\rangle^{\alpha+\frac{2}{\kappa}} \ \cdot\   CK_0 (t'/(m-1)\tilde{T}_{min}(x))^{m-1} \langle x\rangle^{\frac{\alpha}{2}+\frac{1}{\kappa}}  \right.\nonumber\\
 && \left. \qquad\qquad  +  CK_0 (t'/(m-1)\tilde{T}_{min}(x))^{m-1} \langle x\rangle^{\frac{\alpha}{2}+\frac{1}{\kappa}} \ \cdot\ C^4 K_2 \langle x\rangle^{3(\frac{\alpha}{2}+\frac{1}{\kappa})}  \right. \nonumber\\
 && \left. \qquad\qquad 
 \ +\  C^3 K_2^{2/3}  
(t'/(m-1)\tilde{T}_{min}(x))^{\gamma(m-1)/2} \langle x\rangle^{\alpha+\frac{2}{\kappa}} \ \cdot\ C^2 K_1 \langle x\rangle^{\alpha+\frac{2}{\kappa}}   \right\}
   \nonumber\\
&\lesssim &  C^3 K_2^{2/3}  
(t/m\tilde{T}_{min}(x))^{\gamma m/2} \langle x\rangle^{\alpha+\frac{2}{\kappa}} \EEA
for $t\le m\tilde{T}_{min}(x)$, where $\tilde{T}_{min}(x):=\tilde{T}_{min}(0,x)=(C^3
K_2^{2/3} (1+|x|)^{\alpha+\frac{2}{\kappa}})^{-1}$.

\medskip

{\em For larger $t$}, we apply Schauder estimates to eq. (\ref{eq:v}) defining $v^{(m)}$.
Compared to \S 4.2.2, the replacement of $\sup_{Q^{(j)}} |u^{(m)}|$ by $\sup_{Q^{(j)}}
|v^{(m)}|$ leads to an extra prefactor $(t/mT_{min}(t,x))^m\le (t/m\tilde{T}_{min}(t,x))^{\gamma m/2}$. However, due to the right-hand side $f^{(m)-1}(t,x)=-v^{(m-1)}(t,x)\cdot\nabla u^{(m-1)}(t,x)$, there appears an extra contribution in the bound (\ref{eq:S1}) for $|\nabla v^{(m)}(t,x)|$, namely  (concentrating as in \S 4.2.2 on
the main term in the decomposition, for which $2^j\approx \tilde{T}_{min}(t,x)$),

\BEQ 2^{j/2} R_b^{-1}\  \cdot\ 2^{j\gamma/2} ||f^{(m-1)}||_{\gamma,Q^{(j)}} \approx
\tilde{T}_{min}(t,x)^{(1+\gamma)/2} ||v^{(m-1)}\cdot\nabla u^{(m-1)}||_{\gamma, Q^{(j)}}.\EEQ

By induction hypothesis and H\"older interpolation, 
\BEA && ||v^{(m-1)}\cdot\nabla u^{(m-1)}||_{\gamma,Q^{(j)}} \lesssim  ||v^{(m-1)}||_{\gamma,Q^{(j)}} ||\nabla u^{(m-1)}||_{\infty,Q^{(j)}} + ||v^{(m-1)}||_{\infty,Q^{(j)}}
||\nabla u^{(m-1)}||_{\gamma, Q^{(j)}} \nonumber\\
&&  \lesssim  \left( ||v^{(m-1)}||_{\infty,Q^{(j)}}^{1-\gamma} ||\nabla v^{(m-1)}||_{\infty,Q{(j)}}^{\gamma} \right) ||\nabla u^{(m-1)}||_{\infty,Q^{(j)}} + 
||v^{(m-1)}||_{\infty,Q^{(j)}} \left(||\nabla u^{(m-1)}||_{\infty,Q^{(j)}}^{1-\gamma} ||\nabla^2  u^{(m-1)}||_{\infty,Q{(j)}}^{\gamma} \right) \nonumber\\
&&  \lesssim (t/(m-1)\tilde{T}_{min}(t,x))^{\gamma (m-1)/2} \langle x\rangle_t^{(3+\gamma)(\frac{\alpha}{2}+\frac{1}{\kappa})} \left[ (CK_0)^{1-\gamma} (C^3 K_2^{2/3})^{\gamma} C^2 K_1 + CK_0 (C^2 K_1)^{1-\gamma} (C^4 K_2)^{\gamma} \right]
\nonumber\\
&& \le 2 C^{3+2\gamma} (K_2^{2/3})^{1+(1+\gamma)/2} (t/(m-1)\tilde{T}_{min}(t,x))^{\gamma (m-1)/2} \langle x\rangle_t^{(3+\gamma)(\frac{\alpha}{2}+\frac{1}{\kappa})} \nonumber\\
&&\ll C^{3+3(1+\gamma)/2} \left(\langle x\rangle_t^{\alpha+\frac{2}{\kappa}}\right)^{1+(1+\gamma)/2} (K_2^{2/3})^{1+(1+\gamma)/2}  (t/(m-1)\tilde{T}_{min}(t,x))^{\gamma (m-1)/2}  \EEA
 
for $C$ large. Upon multiplication by $\tilde{T}_{min}(t,x)^{(1+\gamma)/2}$
we obtain $ 2^{j/2} R_b^{-1}\  \cdot\ 2^{j\gamma/2} ||f^{(m-1)}||_{\gamma,Q^{(j)}} 
\lesssim C^3 K_2^{2/3}  (t/(m-1)\tilde{T}_{min}(t,x))^{\gamma(m-1)/2} 
\langle x\rangle_t^{\alpha+\frac{2}{\kappa}}$, which is the expected bound for
$|\nabla v^{(m)}(t,x)|$, {\em except} that we still have a factor $(t/(m-1)\tilde{T}_{min}(t,x))^{\gamma(m-1)/2}$ instead of the required $(t/m\tilde{T}_{min}(t,x))^{\gamma m/2}$.

\medskip

\noindent By a minor modification of Proposition \ref{prop:Schauder}, consisting by and large
in substituting $\int_{t'}^t ds\, ||f(s)||_{\infty,Q^{(j)}(s)}$, $t'<t$ to $(t-t')||f||_{\infty,Q^{(j)}}$ where $Q^{(j)}(s)$ is the intersection of the ball $Q^{(j)}$ with
the time-slice $t=s$, in order to take advantage of the extra  factor in $O(1/m)$
coming from the time integral for $s\ll t$, we are able to extract an extra factor
$(t/m\tilde{T}_{min}(t,x))^{\gamma /2}$ for $t\le m\tilde{T}_{min}(t,x)$, where
$\gamma$ is the H\"older exponent. This explains at last why we only obtain 
a prefactor in $O((t/m\tilde{T}_{min}(t,x))^{\gamma m/2})$ in the end for the bound
(\ref{eq:ind-nabla-v}). We do not provide details of this computation since it may
be found in our previous article \cite{Unt-Bur1}, see point (ii) in the proof
of Theorem 3.2.




\section{Appendix}


\begin{Lemma} \label{App:lem:alpha}
Let $A_n$, $n\ge 0$ be a sequence in $\R_+^*$ satisfying an inductive
inequality of the form  $A_{n+1}\le c_1+c_2 A_n^{\alpha}$, with $c_1,c_2>0$ and
$\alpha\in(0,1)$. Then there exists a constant $C_{\alpha}>0$ depending only on $\alpha$ such that $A_n\le \max\left(A_0, C_{\alpha}\max(c_1,c_2^{1/(1-\alpha)}) \right)$ for every $n\ge 1$.
\end{Lemma}

{\bf Proof.} Clearly $A_n\le B_n$, where the sequence $(B_n)_{n\ge 0}$ is defined by
the inductive relation $B_{n+1}=c_1+c_2 B_n^{\alpha}$, with $B_0=A_0$. Let 
$B^*$ be the unique positive fixed point of $\phi:B\mapsto c_1+c_2 B^{\alpha}$. By standard arguments, $(B_n)_{n\ge 1}$ is increasing (resp. decreasing) if $B_1\le B^*$,
resp. $B_1\ge B^*$, and $B_n\to B^*$.  The function $B\mapsto \psi(B):= B-\phi(B)$ 
$(B\ge 0)$
is minimal on $B_*:=(\alpha c_2)^{1/(1-\alpha)}\le c_2^{1/(1-\alpha)}$, 
and increases on the interval $[B_*,+\infty)$. By construction $B^*\ge B_*$ and
$\psi(B^*)=0$.
Let $B_0:=C_{\alpha}\max(c_1,c_2^{1/(1-\alpha)})$. By definition $B_0\in[B_*,+\infty)$. 
We show that $\psi(B_0)\ge 0$, implying $B^*\le B_0$. There are two cases. If
$c_1\ge c_2^{1/(1-\alpha)}$, then $\psi(B_0)\ge (C_{\alpha}-1-C_{\alpha}^{\alpha})c_1$.
In the contrary case, $\psi(B_0)\ge (C_{\alpha}-1-C_{\alpha}^{\alpha})c_2^{1/(1-\alpha)}$. 
Thus in both cases $\psi(B_0)\ge 0$ provided $C_{\alpha}$ is large enough. \hfill\eop



\end{document}